\documentclass{article}
\usepackage[utf8]{inputenc}
\usepackage{amsmath, amssymb, mathtools}
\usepackage{graphicx}
\usepackage{lipsum} 
\usepackage{xcolor}
\usepackage{authblk}
\usepackage{hyperref}
\usepackage{subcaption}
\usepackage{booktabs}
\usepackage{multirow}
\usepackage{epstopdf}
\usepackage{epstopdf}
\usepackage{subcaption}
\usepackage{hyperref}
\usepackage{multirow} 
\usepackage{tabularx}
\usepackage{booktabs} 
\usepackage{siunitx}  
\usepackage{lipsum} 
\usepackage[sort]{natbib}
\usepackage[toc,page]{appendix}
\usepackage[a4paper, total={6.5in, 9in}]{geometry} 

\newcommand{\mean}{\operatorname*{mean}}
\newcommand{\range}{\operatorname*{range}}
\providecommand{\keywords}[1]{\textbf{Keywords:} #1}

 \usepackage{amsmath}

\DeclareMathOperator{\cond}{cond} %

\title{Point Cloud Quality for Meshfree Methods}
\author[1,2]{Mohsen Abdolahzadeh}
\author[1]{Oleg Davydov}
\author[2]{Isabel Michel}
\author[2,3]{Pratik Suchde\thanks{Corresponding author. Email: pratik.suchde@gmail.com}}

\affil[1]{Department of Mathematics, JLU Giessen,
	Arndtstrasse 2, 35392 Giessen, Germany}
\affil[2]{Fraunhofer ITWM, Fraunhofer-Platz 1,
	67663 Kaiserslautern, Germany}
\affil[3]{University of Luxembourg, 2 avenue de l'universit\'e,
	4365 Esch-sur-alzette, Luxembourg}
 \date{\today}

\begin{document}
 	\maketitle

	\begin{abstract}
	Mesh quality is very well studied and widely used to quantify a good mesh.
	In contrast, a systematic study of quality of meshfree point clouds is lacking.
	This gap makes it difficult to substantiate the common claim that
	generating a good-quality point cloud is easier than generating a good mesh.
	Various definitions of point cloud quality have been proposed, some of which 
	have theoretical significance for proving convergence and error bounds,
	while others are used in computational studies.
%
%
%
%
	In this work, we compare and contrast existing point cloud quality metrics and introduce a few new ones.
	We conduct extensive numerical tests with a meshfree collocation method across a wide range of scenarios, including both
	elliptic and hyperbolic equations, 2D and 3D domains, and variations in parameters of the numerical method.
	Based on these tests, we assess which quality metrics best correlate with numerical error.
	Our findings reveal six metrics that consistently serve as reliable indicators of point
	cloud quality, while also demonstrating that several widely used metrics are poor predictors of accuracy.
	\end{abstract}


 	\section{Introduction}
\label{sec:Intro}



Most numerical methods for solving partial differential equations (PDEs) require the creation of a mesh
over the computational domain. Mesh generation, particularly for complex geometries, is widely recognized
as one of the most challenging aspects of mesh based methods, with solution accuracy and stability
heavily dependent on mesh quality.
Meshfree methods have gained increasing attention, in part, to avoid the challenges of mesh generation 
\cite{liu2005introduction}.
These methods, also known as meshless methods, replace mesh generation by the generation of a scattered
set of points, referred to as a \textit{point cloud}, without any topological connectivity.
A common assumption in meshfree literature is that it is easier to generate a good quality point cloud than
a good quality mesh \cite{onate1996finite}.
An overview of different point cloud generation methods can be found in \cite{suchde2023point},
which explores a few aspects of this problem. In this work, we aim to answer a more fundamental question:
\textit{what constitutes a good quality point cloud?}

Just as in mesh-based methods, the quality of the point cloud impacts the accuracy and stability of
the meshfree scheme used. Thus, a clear understanding of what defines a good point cloud is essential.
Poor-quality point clouds can lead to ill-conditioned numerical systems, reduced convergence rates, and
inaccurate solutions. Establishing reliable ways to quantify point cloud quality is thus crucial for
guiding point cloud generation, adaptation, and refinement, ultimately improving the robustness and efficiency
of meshfree methods. This forms the central goal of the present work.

Note that in the context of meshfree methods, the term ``point cloud" as used here, has also been referred to as
particle cloud~\cite{goswami2010interactive},
point set~\cite{liu2009meshfree}, node set~\cite{westermann2025stability},
collocation points~\cite{zhang2001least}, discretization points~\cite{belytschko1999elastic},
data points~\cite{kansa1990multiquadrics} and scattered data~\cite{fasshauer2007meshfree}.

The quantification of mesh quality has been extensively studied over the last few decades, with numerous review
articles and books dedicated to the topic \cite{frey2007mesh, burkhart2013finite,
george2019meshing,lo2002finite,ho1988finite}. Several widely accepted metrics,
such as aspect ratio, orthogonality, and skewness, among others, are commonly used to evaluate mesh quality for both
finite element and finite volume methods. In contrast, there is no detailed study for assessing point cloud quality
in meshfree methods. A number of different metrics have been proposed in meshfree literature to quantify the quality of point
clouds, varying significantly in their underlying principles. Some metrics are introduced
purely for theoretical purposes, such as proving convergence estimates and error bounds, while others are designed
for practical applications, such as improving the conditioning of linear systems. However, the extent to which these
metrics influence numerical accuracy in meshfree methods remains largely unexplored. This work addresses this gap
by providing a comprehensive investigation into the assessment of point cloud quality.

To systematically evaluate point cloud quality, we  examine 25 different quality metrics: 15 from the literature and 
10 introduced in this work.
The existing metrics considered have been introduced for various purposes; some were originally
proposed to assess point cloud quality, while others were developed for point cloud generation algorithms, in particular 
with the idea
of producing ``quasi-uniform'' point distributions \cite{suchde2023point,van2021fast, slak2019generation}.
Here, we adapt them all as point cloud quality indicators. To facilitate the analysis and interpretation of such a large set of metrics,
we introduce a classification system that groups them based on their definitions. Similar to mesh quality assessment,
where metrics are computed per mesh element, all point cloud quality metrics considered operate locally on a neighborhood
of each point.
To quantify global point cloud quality, we apply seven different mappings to aggregate the local metrics
into global ones.

To the best of our knowledge, a systematic study on point cloud quality has not been conducted in
meshfree literature. However, similar work does exist in the context of mesh quality. For instance,
\cite{gao2017evaluating} examines the linear correlation between various mesh quality metrics and
an estimator of the accuracy of the numerical solution of PDE problems. 
This analysis, however, relies solely on Pearson's correlation
coefficient, 
which is limited to measuring linear dependence between the mesh
quality metrics and error estimator.

In contrast to this we compare different metrics using Spearman's rank correlation coefficient $\rho$ that assesses monotonic relationships 
between the metric and the error of the solution of the PDE. Indeed, we see the main application of a quality metric 
in providing an efficient and automatic way to decide whether the point clouds' quality increases after any modifications
aimed at improving it. Assuming that the ground truth quality is given by the size of the error of the numerical solution,
a quality metric fulfills this purpose if it is comonotonic with the error, making the monotonic correlation measures such as 
Spearman's $\rho$ or Kendall's Tau appropriate tools. Measuring the linear correlation between the error and a
metric would only be meaningful if the metric is appropriately normalized, which may be possible for some metrics, but not in
general since complicated nonlinear dependencies  are common, as our experiments in Section~\ref{sec:Validation_Results} indicate.

We restrict our study to PDEs with known
analytical solutions, allowing for precise error quantification. This allows us to systematically
evaluate the correlation between numerical error and point cloud quality metrics.
To assess the robustness of different metrics, we conduct extensive numerical experiments
covering elliptic and hyperbolic PDEs, domains in $\mathbb{R}^2$ and $\mathbb{R}^3$, and
variations in numerical parameters such as resolution, support size, and polynomial order.
By analyzing these diverse cases, we identify several quality metrics that consistently serve as
reliable indicators of numerical accuracy.


Since theoretical analysis of the relationship between point cloud quality and solution error is
challenging in meshfree methods, we adopt the numerical approach described above. We believe that this
serves as a crucial first step toward understanding the impact of point cloud quality on solution accuracy.
Our work focuses on meshfree collocation methods, specifically we consider two types of the meshfree generalized
finite difference methods (GFDM), which is a strong-form meshfree approach based on weighted least squares approximations.
This class of methods has been widely applied to various scientific and engineering problems (\cite{veltmaat2022mesh,
suchde2018meshfree, michel2021meshfree, bharadwaj2022discrete, leithauser2024predicting}).
Our systematic study also provides a framework that can be readily extended to other meshfree methods.

The paper is organized as follows.
Section~\ref{metrics} introduces the point cloud quality metrics considered in this study and discusses
their classification and aggregation maps that generate global metrics.
Section~\ref{methodology} describes our methodology for assessing the performance of quality metrics,
including generation of  point clouds for testing and methods to estimate  
the monotonic correlation  between quality metrics and numerical error.
Section~\ref{sec:PDEs_GFDM} presents the types of PDE  considered and the numerical approximation methods employed
in this study.
The validation of our methodology is presented in Section~\ref{sec:Validation_Results}, including a comparison
to the mesh-based case.
Section~\ref{eval} then presents our numerical experiments to evaluate the various quality metrics.
In particular, a comprehensive  discussion of the results is given in Section~\ref{sec:Discussion}.
Finally, Section~\ref{sec:Conclusion} concludes the paper and outlines directions for future work.
For the convenience of the reader some technical information, such as detailed definitions of the metrics, 
 is pushed to the Appendix.
%


 	
\section{Metrics to assess point cloud quality}\label{metrics}

Consider a computational domain $\Omega\subset \mathbb{R}^\nu$ with boundary $\partial \Omega$.
The domain $\overline{\Omega} = \Omega \cup \partial \Omega$ is discretized with
a finite set of $n$ points, \(\mathbf{X}:= \{\mathbf{x}_i\}_{i=1}^{n} \), referred to as a \emph{point cloud}.
At a point $\mathbf{x}_i \in \mathbf{X}$, all approximations are carried out on a \emph{neighborhood} or
\emph{support} \(\mathbf{S}_{i} \subset \mathbf{X}\) of a small number of nearby points that includes $\mathbf{x}_i$.
This system of supports is the only topological structure imposed on the point cloud.

\medskip

\noindent
\textbf{Remark } Meshfree literature includes a wide range of terminology.
Points $\mathbf{x}_i$ are also often referred to as \emph{nodes} or \emph{particles}.
The support  \(\mathbf{S}_{i}\) of $\mathbf{x}_i$ is also known as a \emph{support domain} or a \emph{set of influence}.
The point cloud \(\mathbf{X}\) is also referred to under various different terms, as discussed in Introduction.

\subsection{Local and global quality metrics}\label{locaAndGlobalDeff}


In mesh-based methods, mesh quality is typically assessed by local metrics on individual elements,
such as the aspect ratio of a triangular element.
A quantification of the quality of the entire mesh is obtained through an aggregated metric, such as the maximum, minimum or an
average of these local values.
This facilitates local mesh improvement algorithms that adjust poor-quality regions \cite{knupp2020fundamentals, knupp2022worst},
while also providing an overall assessment of the mesh quality.

We aim to use a similar approach in the meshfree setting.
For a point \(\mathbf{x}_i \in \mathbf{X}\),
we consider a \emph{local quality metric} $d_i$ defined on the support $\mathbf{S}_{i}$.
The \emph{global quality metric} \(d_{ \mathbf{X}}\) that quantifies the quality of the point cloud $\mathbf{X}$
as a whole is obtained by applying an \emph{aggregation map} $\mu$
to the sequence of local metrics $\{d_i:\mathbf{x}_i \in \mathbf{X} \}$, or an appropriate subsequence of it, for example
excluding boundary points whenever Dirichlet boundary conditions are used, see Section~\ref{sec:Procedure}.
%
In this work, we consider seven aggregation maps to obtain the global quality metrics:
maximum, minimum, arithmetic mean, standard deviation, root mean square, median  and range,
denoted as \texttt{max}, \texttt{min}, \texttt{mean},  \texttt{std}, \texttt{rms},
\texttt{median} and \texttt{range}, respectively.
We use standard definitions of these maps, with the range defined as the difference between
maximum and minimum values.

%

\subsection{Types of local quality metrics}\label{typesAndMetrics}

Various criteria have been proposed in meshfree literature to assess point cloud quality.
These include, in particular, geometric properties such as the minimum distance between points,
physical characteristics of the point cloud considered as a system of interacting particles, or
characteristics of a mesh built by tessellating the point cloud.

In this work, we examine 25 quality metrics to assess local point cloud quality, including
existing metrics from the literature, metrics repurposed from point
cloud generation methods or from a theoretical error analysis, and a few new metrics introduced here.
A full list of metrics is provided in Table~\ref{tab:qualityMetrics}. We emphasize that these metrics are computed
locally within the support $\mathbf{S}_{i}$ of each point
$\mathbf{x}_i$, even if some of them are usually defined for the entire point cloud $\mathbf{X}$ in the literature.
Most local metrics considered depend on all points in \(\mathbf{S}_i\), while some use only a 
subset of \(\mathbf{S}_i\) consisting of the closest points. In particular, the 
mesh-based metrics (see the classification introduced below) depend only on such a subset, whereas 
all other metric types use the entire support.

To make it easier to navigate the many metrics considered here and to better present
their differences, we classify local quality metrics into several groups.
This classification is based on how each metric is
defined and which attributes of the local point cloud it considers.
%
%
Based on these considerations, we introduce five \emph{types} of point cloud quality metrics

\begin{itemize}
	\item \textbf{Energy-based metrics}: By modeling a point cloud as an interacting particle system, these
	metrics assess point cloud quality via the system's energy defined either through physical principles
	or in an \emph{ad hoc} manner. A better point cloud is expected to have lower energy.

	\item \textbf{Equilibrium-based metrics}: 
	These metrics are defined by evaluating the imbalance of forces in a system, based on an energy function.
	They are computed using the norm of the gradient of the energy function,
	with lower values indicating a more stable configuration.

	\item \textbf{Geometry-based metrics}: These metrics evaluate point cloud quality based directly on
	geometric properties such as nearest and farthest neighbors, pairwise distances,
 	and similar attributes.

	\item \textbf{Mesh-based metrics}: These metrics rely on characteristics of mesh  elements associated
	with a point cloud. They require meshing of the local supports $\mathbf{S}_{i}$ unless a global
	mesh of the point cloud $\mathbf{X}$ is available. 
%
	\item \textbf{System matrix-based metrics}:
	These metrics assess point cloud quality based on the local characteristics of the system matrix
    associated with the numerical scheme. For collocation methods, this is constructed from numerical
	differentiation weights that approximate derivatives at each point (see Section~\ref{numericalMethods}).
\end{itemize}


For the sake of brevity, we only list all the metrics here,
and defer their descriptions to Appendix~\ref{qualityMetricDetails}.
Table~\ref{tab:qualityMetrics} lists all the quality metrics considered in the present study,
and provides reference to their descriptions in Appendix~\ref{qualityMetricDetails}.
The table also provides information on each metric's source,
including the work where they were used, and for what purpose.
Metrics marked with ``present study'' indicate metrics introduced in this work.
The table also provides an abbreviation for each metric, which will be used throughout this work.


\begin{table}[!htp]
	\caption{Summary of the point cloud quality metrics used.
	\emph{Type} refers to the classification of metrics outlined in Section~\ref{typesAndMetrics},
	\emph{Description} gives a reference to the section of  Appendix~\ref{qualityMetricDetails} where
	the metric is explained in detail, \emph{Origin} cites a
	literature source where the metric has been used previously, although we have not investigated in all cases
	where the respective methods were originally proposed,
	with \emph{present study} indicating metrics introduced here.
	\emph{Source} outlines the original  application of each metric adopted from the literature.
	The \emph{Abbreviations}  of the last column are used everywhere in the paper to refer to particular metrics.
	}
	\centering
	\resizebox{\textwidth}{!}{%
		\begin{tabular}{lcccccc}
			\toprule
			Type & Metric & Description &Origin &
			Source &
			Abbreviation \\
			\midrule
			\multirow{5}{*}{Energy}
			& Riesz Kernel Energy &   \ref{rieszEnergy}       &
			\cite{vlasiuk2018fast}        &
			Point cloud generation              & RKE \\
			& Log-Kernel Energy  &   \ref{logaEnergy}        &
			\cite{kunc2019generation}        &
			Point cloud generation              & LKE \\
			& Inter-Molecular Energy &   \ref{IMenergy}      &
			\cite{zhang2005node}   &
			Point cloud generation              & IME \\
			& Wendland Kernel Energy &   \ref{WandlandEnergy}    & Present
			study                                      &
			-                 & WKE \\
			& Centroid Energy &   \ref{balaceKernelEnergy}         & Present
			study                                      &
			-                 & CE  \\
			\midrule
			\multirow{5}{*}{Equilibrium}
			& Riesz Kernel Potential Force&   \ref{rieszPF} & Present
			study                                    &
			-                 & RKPF\\
			& Log-Kernel Potential Force&   \ref{logPF}  & Present
			study                                    &
			-                 & LKPF\\
			& Inter-Molecular Potential Force&   \ref{potentialForceIMP} &
			\cite{negi2021algorithms} &
			Point cloud generation              & IMPF\\
			& Wendland Kernel Potential Force &   \ref{wandPF} & Present
			study                                  &
			-                 & WKPF\\
			& Centroid Potential Force &   \ref{baPF}             & Present
			study                                      &
			-                 & CPF  \\
			\midrule
			\multirow{6}{*}{Geometry}
			& Separation Distance &   \ref{Separation}    &
			\cite{fornberg2015fast}       &
			Quality metric                 & SD  \\
			& Fill Distance   &   \ref{fillDistance}            &
			\cite{van2021fast, du2002meshfree}          & Point cloud generation, Quality metric                 & FD
			\\
			& Uniformity Ratio  &   \ref{meshRatio}               &
			\cite{van2021fast}     & Quality metric                 & UR  \\
			& Packing Density   &   \ref{SecPackDensity}          &
			\cite{van2021fast}     &
			Quality metric, Point cloud generation & PD \\
			& Average of Local Regularity &   \ref{localRegularity}    &
			\cite{van2021fast}     &
			Quality metric                 & ALR \\
			& Range of Local Regularity &   \ref{localRegularity}    &
			\cite{van2021fast}     &
			Quality metric                 & RLR \\
			\midrule
			\multirow{3}{*}{Mesh}
			& Inverse Aspect Ratio    &  \ref{aspectRationAR}            &
			\cite{geuzaine2009gmsh}      & Mesh
			quality metric            & IAR  \\
			& Determinant of the Jacobian Matrix & \ref{jacobM}                  & Present
			study                                      & -                 & DJM  \\
			& Mass Centers Deviation &   \ref{massCenterD}    &
			\cite{fu2019optimal}                    &
			Point cloud generation              & MCD \\
			\midrule
			\multirow{6}{*}{System Matrix}
			& Growth Function via $\ell_1$-norm   &\ref{growthF}          & \cite{davydov2018minimal}
			                                      &
			Error bounds                 & GFL1  \\
			& Growth Function via $\ell_2$-norm   &\ref{growthF}          & \cite{davydov2018minimal}
			                                      &
			Error bounds                 & GFL2  \\
			& Condition Number&\ref{conditionN}           & Present
			study                                      &
			-                 & CN  \\
			& Recovery
			Error &\ref{maxAppErro}            & \cite{davydov2021approximation}
		                                      & Error bounds                 & RE  \\
			 & Stability  &\ref{stabilityMe}       & Present
			study                                      &
			-                 & ST   \\
		 & Diagonal Weights  &\ref{sysMat}& Present
			study                                &
			-                 & DW  \\
			\bottomrule
		\end{tabular}%
	}
	\label{tab:qualityMetrics}
\end{table}

We note that almost all metrics considered here are non-negative, with $d_i\ge 0$.
The only metric that accepts negative values is the ``Diagonal Weights'' metric,
see Appendix~\ref{qualityMetricDetails} for more details.

Each of the $7$ aggregation map functions defined in Section~\ref{locaAndGlobalDeff} is applied
to all $25$ local metrics mentioned in Table~\ref{tab:qualityMetrics}. This results
in a total of \(25 \times 7 = 175\) global metrics considered in this study.

 	\newpage

\section{Methodology}\label{methodology}

In this section, we introduce our methodology for evaluating the performance of the quality metrics
described in Section~\ref{metrics}. Our goal is to identify which metrics best predict the performance of
meshfree collocation methods for solving PDEs.
Clearly, robust and reliable quality metrics are essential for designing effective algorithms to generate or improve
point clouds.

In order to compare performance of different metrics, we choose several PDEs with known analytical solutions.
We generate a large number of point clouds of varying quality and compute the numerical solution using a
meshfree collocation method. An ideal quality metric would always assign a higher quality score to point clouds
that produce lower errors of the numerical solution against the analytical solution.
However, realistically, no quality metric may achieve such an ideal performance,
as the numerical error depends not only on point cloud quality, but also on several additional factors,
such as variations in the smoothness of the PDE solution across the computational domain. Nevertheless, we compare different
metrics on the degree to which they predict the numerical error, in order to  identify which metrics achieve this goal most
reliably and consistently.


\subsection{Point cloud generation}\label{nodSetGeneration}
As discussed earlier, we evaluate the robustness and reliability of quality metrics described
in Section~\ref{metrics} by assessing their performance across a large number of point clouds.
A key consideration here is that point clouds must have the same number of points
to ensure a meaningful comparison.

\begin{figure}[!ht]
	\centering
	\begin{subfigure}[!htp]{0.45\textwidth}
		\centering
		\includegraphics[width=\textwidth]{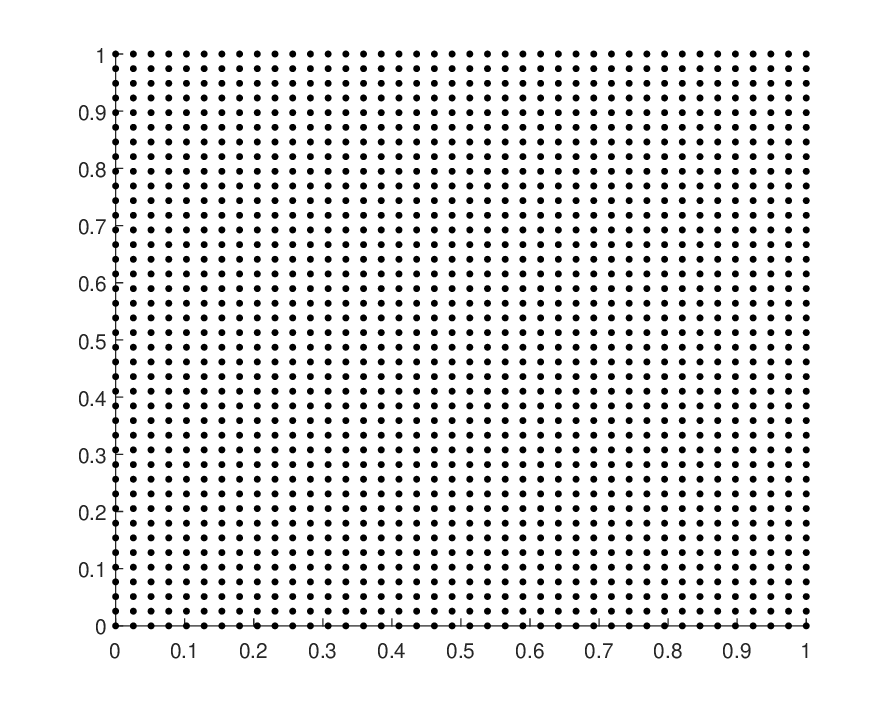}
		\caption{Unperturbed Cartesian point cloud, \(\alpha = 0\)}
	\end{subfigure}
	\begin{subfigure}[!htp]{0.45\textwidth}
		\centering
		\includegraphics[width=\textwidth]{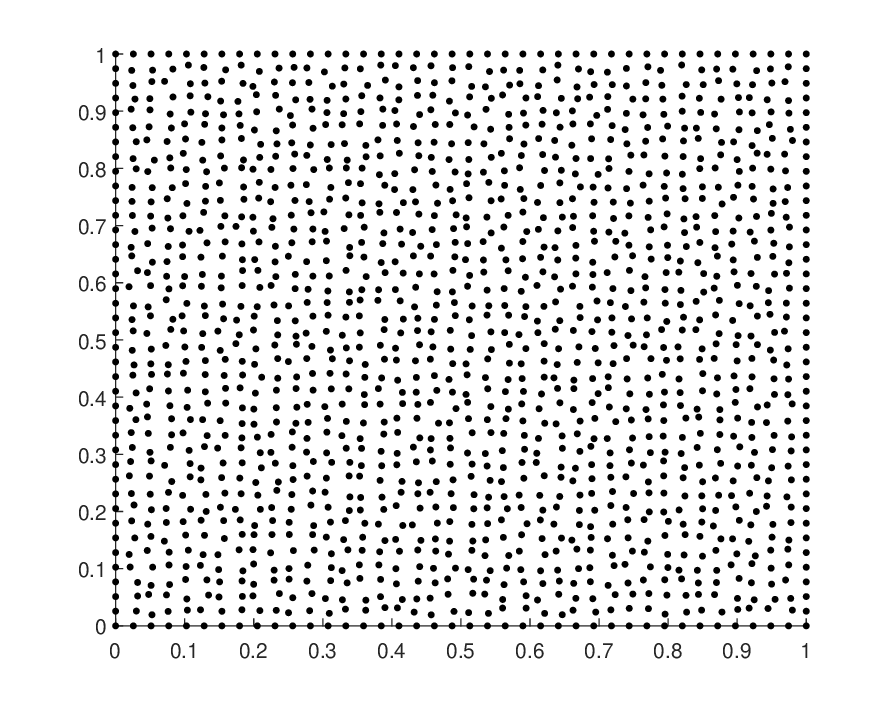}
		\caption{\(\alpha = 0.25\)}
	\end{subfigure}
	\\
		\begin{subfigure}[!htp]{0.45\textwidth}
		\centering
		\includegraphics[width=\textwidth]{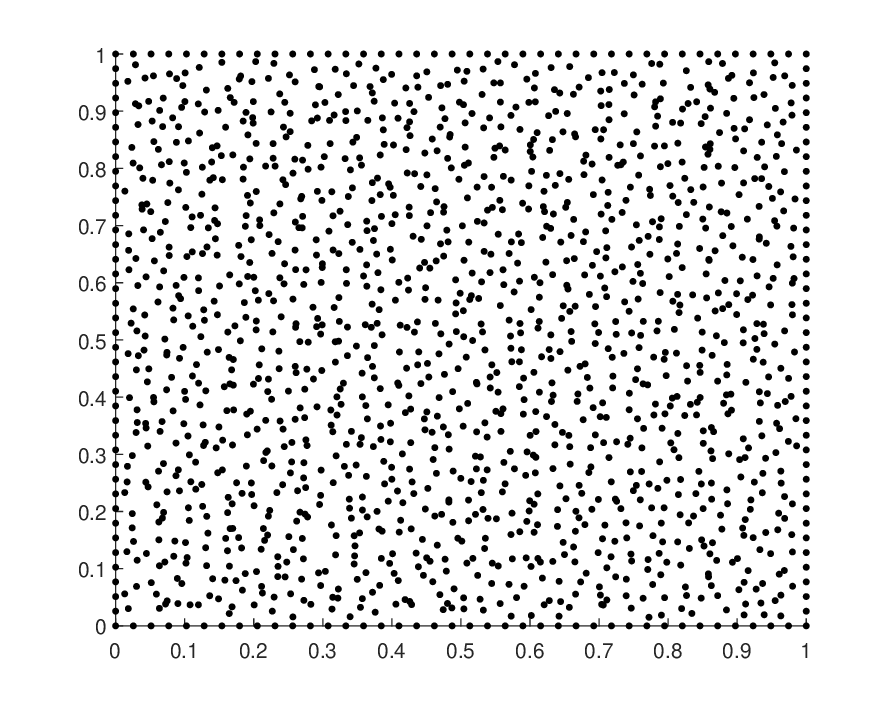}
		\caption{\(\alpha = 0.50\)}
	\end{subfigure}
		\begin{subfigure}[!htp]{0.45\textwidth}
		\centering
		\includegraphics[width=\textwidth]{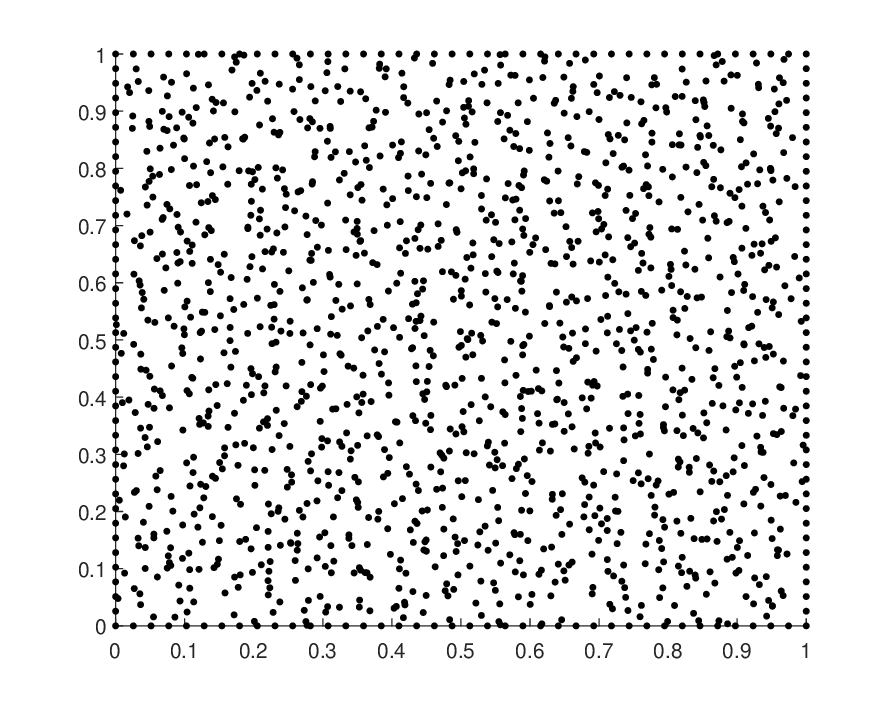}
		\caption{\(\alpha = 1.00\).}
	\end{subfigure}

	\caption{Illustration of point clouds generated for a unit square domain in $\mathbb{R}^2$.}
	\label{visualPointSet}
\end{figure}

We start with a Cartesian grid of $n$ points, with uniform spacing \(h\), that discretizes the
computational domain, denoting it by
\(\mathbf{X}^0 := \{\mathbf{x}^0_i\}_{i=1}^{n} \subset \mathbb{R}^\nu\).
Each point in this initial point cloud is randomly perturbed, 
which generates a new point cloud
\(\mathbf{X}^{\alpha} := \{\mathbf{x}^{\alpha}_i\}_{i=1}^{n}\)
created as
\begin{equation}
	\mathbf{x}^{\alpha}_i = \mathbf{x}^0_i + \alpha \, h \, \boldsymbol{\delta}_i, \quad
	i = 1, \dots, n, \label{setGeneration}
\end{equation}
where 
$\boldsymbol{\delta}_i = (\delta_{i,1}, \delta_{i,2}, \dots, \delta_{i,\nu})$
is a $\nu$-dimensional random variable, where each component is an independent,
uniformly generated
random number $\delta_{i,j} \sim \text{Uniform}(-1, 1)$ in the interval \([-1, 1]\).
%
Furthermore, the factor \(\alpha\) bounds the magnitude of the perturbation.
We generate multiple point clouds \(\mathbf{X}^{\alpha} \) with different values
of \(\alpha\) to obtain point clouds with diverse point configurations. For this,
\(\alpha\) is varied as
\begin{equation}
	\alpha = \alpha(k) = 0.05k, \quad k=1,\ldots,20.
	\label{randomness}
\end{equation}
Thus, for example $k=1$ generates a point cloud with a maximum $5\%$ deviation from the Cartesian
grid,  $k=10$  means a maximum $50\%$ deviation from the Cartesian grid,
and so on. Figure~\ref{visualPointSet} illustrates four such point clouds, with different amounts
of perturbation,  \(\alpha =0,  0.25, 0.5 \text{ and } 1.0\) of the nodes of a $40 \times 40$
Cartesian grid. The figure demonstrates that
the irregularity of the point cloud becomes more pronounced as \(\alpha\) increases.

We note that $\alpha<0.5$ for \( k < 10 \), and the minimum distance \( (1 - 2\alpha)h \) is
ensured between each pair of points. However, for \( k \geq 10 \) we enforce a lower bound on the minimum distance,
ensuring a separation of \(10^{-5} h \), in order to avoid (almost) coincident points.
If two points come closer than this threshold during the generation of a perturbed point cloud,
we regenerate the random perturbation of one of them until the minimum distance is satisfied.

To account for the variability in the random perturbation process we generate multiple point clouds (namely, 50) for each
level of perturbation $\alpha(k)$ in Equation~\eqref{randomness}.
This helps to reduce the influence of the outliers and ensures a more reliable assessment of the
quality metrics. This also ensures that we are
evaluating the performance of the quality metrics across a diverse set of point configurations.
Thus, for each test problem we generate a total of $1001$ point clouds including the unperturbed Cartesian grid.


\subsection{Assessing the predictive power of quality metrics}
\label{sec:monotonicity}

Having generated a diverse set of point clouds as outlined in the previous subsection,
our goal is to determine how well different quality metrics correlate with the numerical error.
As mentioned above, no
metric can perfectly predict numerical error due to the presence of other factors
influencing numerical accuracy. However, some metrics may still serve as
reliable indicators thanks to a consistently high correlation with the error trends.

To objectively assess the correlation between numerical error and a quality metric
across all the point clouds considered, a quantitative approach is essential.
Simple visual inspection of the relationship is insufficient, as it may overlook subtle
trends or inconsistencies. To achieve this, we quantify the monotonic relationship
between each global metric $d_{\mathbf{X}}$ and the numerical error \(e_{\mathbf{X}}\) over all 1001
point clouds $\mathbf{X}$ included in the tests. A strong monotonic (or counter-monotonic) dependence indicates that the
metric is a good predictor of point cloud
quality, whereas a weak dependence suggests that it is less reliable. As discussed in Introduction, monotonic correlation
is an appropriate criterion in contrast to a linear correlation that cannot reveal nonlinear dependencies
between a metric and the error.

To evaluate the monotonic relationship, we consider five different measures of monotonicity
and first evaluate their performance for our setting, see Section~\ref{sec:Validation_Results}.
Specifically, we examine the
%
%
\begin{itemize}
\setlength{\itemsep}{1pt}
\setlength{\parskip}{0pt}
\setlength{\parsep}{0pt}
	\item Spearman's rank correlation coefficient $\rho$ \cite{Gibbons1985},
	\item Kendall's Tau coefficient $\tau$ \cite{Kendall1990},
	\item Minimum Covariance Determinant (MCD) measure \cite{tabatabai2021introduction},
	\item Monotonicity Coefficient (MC) \cite{kachapova2012measure}.
	\item In addition, we include  Pearson's correlation coefficient $r$ \cite{stigler1989francis} that is not specifically designed to test
	monotonicity, but has been used in a relevant study for mesh quality metrics \cite{gao2017evaluating}.
\end{itemize}

To determine which of these monotonicity measures is best suited for our purpose, we compare
in Section~\ref{sec:Validation_Results} their
predictions against observed trends in scatter plots of error versus metric values.
Specifically, we assume that the monotonicity measures that contradict clear visual trends
are not well suited to determine the relation between quality metrics and numerical error.
We observe empirically  that the Kendall's Tau \( \tau \) and Spearman's rank correlation coefficient
\( \rho \) consistently reflect the expected relationships, making them well-suited for our
analysis. In contrast, the remaining three measures exhibited poor agreement with visual
trends. 
We refer to the books \cite{Gibbons1985,Kendall1990}  for a detailed description of
Spearman's rank correlation and Kendall's Tau.

The main tests in Sections~\ref{parameterStudy}--\ref{sec:HyperbolicResults} evaluate the performance
of the quality metrics by Spearman's rank correlation coefficient \(\rho\), even if we could instead use 
the Kendall's Tau that delivers very similar results.

Note that we only consider the absolute values of all monotonicity measures throughout this work, and
whenever an average of \(\rho\) is used, it always refers to the average of the absolute values.
In statistical literature, monotonicity measures are typically
taken as signed values, where positive values indicate co-monotonicity (error increases with increasing metric value)
and negative values indicate counter-monotonicity (error decreases with increasing metric value). However,
since we are concerned only with the strength of the monotonic correlation, whether positive or negative, we report
only absolute values when assessing the reliability of point cloud quality metrics.

	\section{Test problems and computational setup}
\label{sec:PDEs_GFDM}

We evaluate the effectiveness of the quality metrics presented in Section~\ref{metrics}
in two settings of an elliptic and a hyperbolic PDE.
These two types of equations present distinct numerical challenges,
allowing us to assess how quality metrics perform across different numerical scenarios.

\subsection{Elliptic PDE: Poisson equation}
\label{sec:Poisson}
For the elliptic setting, we consider the Poisson equation with Dirichlet boundary conditions
\begin{align}
    \begin{split}
    	\Delta	u &= f  \quad  \text{in} \quad  \Omega, \\
    	u &= g \quad \text{on} \quad \partial \Omega,
    	\label{PoissonEq}
    \end{split}
\end{align}
where \(\Delta\) represents the Laplace operator. In the computational tests in the next sections,
we consider the computational domain \(\Omega\) to be a unit square $[0,1]^2$ in $\mathbb{R}^2$,
and a unit cube $[0,1]^3$ in $\mathbb{R}^3$.
To cover a range of cases, we consider two analytical solutions in $\mathbb{R}^2$,
and two for the $\mathbb{R}^3$ domain. The analytical solutions and corresponding right-hand side
functions $f$ are summarized in Table~\ref{testsFunction}. Dirichlet boundary
conditions $g$ are obtained from the analytical solutions.
\begin{table}[!htp]
\centering
	\caption{Poisson problems: Analytical solutions considered and the corresponding right-hand sides
	for the $\mathbb{R}^2$ and $\mathbb{R}^3$ domains (unit square and unit cube, respectively).}
\begin{tabular}{@{}ccc@{}}
		\toprule
	Dimension	&Exact Solution (\(u\)) & Right-Hand Side (\(f\)) \\
	\midrule
		\multirow{2}{*}{2D}
	&	\(u_{1} = e^{xy}\) & \(f_{1} = (x^2 + y^2)e^{xy}\) \\
	&	\(u_{2} = x^3 \sin(y) + y \cos(x)\) & \(f_{2} = 6x \sin(y) - y \cos(x) - x^3 \sin(y)\) \\
	\midrule
		\multirow{2}{*}{3D}
	&	\(u_{3} = \sin(x)\sinh(y)\sin(z)\) & \(f_{3} = -\sin(x)\sinh(y)\sin(z)\) \\
	&	\(u_{4} = x^3 \sin(y) + y \cos(x) + y \sin(z)\) & \(f_{4} = 6x \sin(y) - y \cos(x) - x^3
	\sin(y) - y \sin(z)\) \\
	\bottomrule
	\end{tabular}
	\label{testsFunction}
\end{table}

To evaluate the performance of different metrics for the Poisson equation,
we use the relative maximum error defined as
\begin{equation}
	e_{\infty} = \frac{\|\hat{u} - u\|_{\infty}}{\|u\|_{\infty}},
	\label{errirInf}
\end{equation}
where $u$ is the analytical solution, $\hat{u}$ is the numerical solution, and the $\infty$-norm is computed as the maximum of the absolute value of a function on
the point cloud $\mathbf{X}$.
%
\subsection{Hyperbolic PDE: Advection equation}
\label{sec:Advection}
For the hyperbolic setting, we consider an advection equation with Dirichlet boundary conditions
\begin{equation}
	\begin{split}
		\frac{\partial \phi}{\partial t} + \mathbf{V} \cdot \nabla  \phi &= 0,
		\quad  \text{in} \quad  \Omega, \\
		\phi &= g \quad \text{on} \quad \partial \Omega, \\
		\phi_0 &= f \quad \text{in} \quad \Omega.
	\end{split}
	\label{hyperbolicPDE}
\end{equation}
where \( \phi = \phi(\mathbf{x}, t) \) is the scalar quantity being transported,
\( \phi_0 := \phi(\mathbf{x}, 0) =f(\mathbf{x})\) is the initial condition,
\( \mathbf{V} = \mathbf{V}(\mathbf{x}) \) is the prescribed velocity field, and
\( \nabla = \left( \frac{\partial}{\partial x}, \frac{\partial}{\partial y} \right) \)
is the gradient operator. The shifted unit square $[-0.5,0.5]^2$ in $\mathbb{R}^2$
is the only computational domain considered here.

We consider a rotational velocity field, \( \mathbf{V} = (y, -x) \), with the exact solution $\phi$
representing a rotation of a Gaussian bell.
The Dirichlet boundary condition \( g \) and the initial condition \( f \) are chosen
such that the exact solution is:
%
%
\begin{equation}
	\phi(x,y,t) =
	\begin{cases}
		\dfrac{\exp \left( - (x - \frac{\cos t}{4})^2 - (y -  \frac{\sin t}{4})^2 \right)
			- \exp (-\frac9{200})}
		{1 - \exp (-\frac9{200})},
		& \text{if } (x - \frac{\cos t}{4})^2 + (y - \frac{\sin t}{4})^2 \leq \frac9{200}, \\[10pt]
		0, & \text{otherwise}.
	\end{cases}
	\label{hyperbolic}
\end{equation}

To evaluate the performance of different quality metrics for the advection equation,
we use the relative maximum error against the analytical solution as in Equation~\eqref{errirInf} at time $t = 3.14$. 
%


\subsection{Meshfree finite differences}\label{numericalMethods}

Throughout this work, we consider the class of generalized finite difference methods (GFDMs) going back to~\cite{jensen1972finite},
where the discretization of a PDE is
obtained by a meshfree collocation with local polynomials. See~\cite{suchde2018conservation} for details and further references.
Consider a point \( \mathbf{x}_i\) with support \(  \mathbf{S}_i \). A general linear
differential operator \( D \) at this point is approximated by a numerical differentiation formula as
\begin{equation}
	D f(\mathbf{x}_i) \approx \sum_{\mathbf{x}_{j} \in \mathbf{S}_i } w_{i,j}
	f(\mathbf{x}_{j}).
	\label{diffFormula}
\end{equation}
%
The weights $w_{i,j}$ are usually determined numerically by solving a linear system, or
minimizing a quadratic functional
of these weights, subject to Equation~\eqref{diffFormula} retaining exactness for polynomials
of a specified order $p$, that is, of degree $p-1$. In all cases the implementation involves
solving a linear system of a small size depending on the number of points in $\mathbf{S}_i$. 
Several choices of the minimization functional
have been considered in the GFDM literature. Note that the same weights are also obtained by applying differential operators  
to a polynomial weighted least squares fit to the data at $\mathbf{S}_i$, which can be related to a 
Taylor series expansion as well, see for example \cite{suchde2018conservation,davydov2018minimal}.

We use two different weighted least squares norms as minimization functionals. For the Poisson problem, we minimize
\(\lVert \mathbf{w}_i \lVert_{2,\mu}^2\ := \sum_{\mathbf{x}_{j} \in \mathbf{S}_i}
w_{i,j}^{2} \rVert \mathbf{x}_{j} - \mathbf{x}_i \lVert^{2\mu}\), where $\mathbf{w}_i$ denotes the vector with the components
$w_{i,j}$ for all $j$ such that $\mathbf{x}_{j} \in \mathbf{S}_i$.
The parameter $\mu$ is chosen equal to the polynomial order, $\mu = p$,
as suggested in \cite{davydov2018minimal} based on the estimates of the numerical differentiation error. 
For this minimization, we compute the weight vectors $\mathbf{w}_i$ using the
open-source code \texttt{mFDlab} \cite{mFDlab}.

For the advection problem, we use the commercial software MESHFREE \cite{meshfree}
developed by the Fraunhofer ITWM and SCAI institutes,
where the differentiation weights are obtained by the minimization of the least squares functional
$\sum_{x_{j} \in \mathbf{S}_i} \frac{w_{i,j}^{2}}{\kappa_{i,j}^{2}} $ defined through
a Gaussian kernel function
$\kappa$ \cite{suchde2018conservation},
\begin{equation}
\kappa_{i,j} =	\exp\left(-\mu \frac{ \|\mathbf{x}_i - \mathbf{x}_j\|^2 }{ \lambda^2 } \right),
\end{equation}
where $\mu $ is a positive constant and $\lambda$ is the support size, see Equation \eqref{range}.
We use $\mu = 4.0$ throughout this work.
Furthermore, to handle the advection term, we use a MUSCL reconstruction scheme with a
SUPERBEE limiter. Time integration is done using the SDIRK2 method \cite{alexander1977diagonally}, which is a
second order implicit time integration scheme. For more details, we refer to \cite{seifarth2018numerische}.


%
%

Once the differentiation weights have been computed, the PDE is discretized using the discrete
derivatives in the right hand side of \eqref{diffFormula}. For example, for the Poisson Equation~\eqref{PoissonEq}, the discrete system to
compute the numerical solution $\hat{u}$ is formed as
\begin{align}
\begin{split}
	\sum_{\mathbf{x}_j \in \mathbf{S}_i} w_{i,j} \hat{u}_{j} & = f(\mathbf{x}_i),
	\qquad  \mathbf{x}_i \in \Omega,
	\\
	u(\mathbf{x}_i) & = g(\mathbf{x}_i), \qquad  \mathbf{x}_i \in \partial \Omega,
		\label{discritPoisson}
\end{split}
\end{align}
where the weights $w_{i,j}$ are from Equation~\eqref{diffFormula} for $D=\Delta$.
A similar procedure is used for the spatial discretization of the advection equation, using the discrete gradient operators.

\subsection{Numerical procedure}
\label{sec:Procedure}

For both the Poisson and advection problems, we evaluate the effectiveness of all the metrics
mentioned in Section~\ref{metrics} using the procedure outlined in Section~\ref{methodology}.
We follow the same procedure for the 2D advection problem and both 2D and 3D Poisson problems,
ensuring a fair comparison of the quality metrics across different numerical settings.
For each domain, we generate a total of $1001$ point clouds by randomly perturbing
a Cartesian grid, as described in Section~\ref{nodSetGeneration}.
In particular, the same set of point
clouds is used for both 2D cases, with the only difference being that point clouds on $[0,1]^2$
for the Poisson problem are shifted to $[-0.5,0.5]^2$ for the advection problem.
The effectiveness of each point cloud quality metric  is then determined using the monotonicity
measures introduced in Section~\ref{sec:monotonicity} to estimate their  correlation with the numerical error metric values.

Since Dirichlet boundary conditions are imposed for both PDEs, random perturbations are applied
only to interior points, while boundary points remain fixed, as illustrated by
Figure~\ref{visualPointSet}. Furthermore, all aggregation maps to compute the global quality metrics
are computed over local metrics on interior points only. For example, when the global metric is defined as
the mean of the local metric values at all points, the mean is taken only over non-boundary points.

 	\section{Validation of methodology}
\label{sec:Validation_Results}

In this section, we validate our approach to assessing point cloud quality metrics
using the Poisson problem in $\mathbb{R}^2$. We first check how our approach presented in Section~\ref{methodology} applies
to a standard mesh quality metric: the aspect ratio in the finite element method, and then compare the performance of the five
monotonicity measures of Section~\ref{sec:monotonicity} in the meshfree setting.

\subsubsection*{Mesh quality in Finite Element Method}

Mesh quality has been widely studied and is well understood, with well-established
metrics used to assess the impact of the mesh on numerical accuracy.
This offers a straightforward way to validate our methodology for point cloud quality
assessment by first applying it to mesh quality. Specifically, we use our approach
outlined in Section~\ref{methodology}
to analyze one of the most commonly used mesh quality metrics, the aspect ratio.
For a triangle $T$ in a 2D triangular mesh, the inverse aspect ratio \(\gamma_{T}\) of the element $T$ is given by
\begin{equation}
	\gamma_{T} = \frac{2\sqrt{3}\rho_T}{h_T},
   \label{aspectRatio}
\end{equation}
where $h_T$ represents the diameter of $T$ (length of the longest edge), and $\rho_T$ denotes the inradius \cite{george2019meshing}. 
We refer to \(\gamma_{T}\) as the \emph{inverse aspect ratio} because it is the reciprocal of the more common
aspect ratio convention \cite{geuzaine2009gmsh}. 
For more details, see Appendix~\ref{aspectRationAR}.

The global aspect ratio as a mesh-quality metric is determined as the minimum over all elements, $\text{min}(\gamma_{T})$, 
which is the standard aggregation map for this metric in the mesh generation
literature.

To validate our methodology, we consider the Poisson problem \eqref{PoissonEq}
on the unit square in \(\mathbb{R}^2\) and exact solution \(u_{1}\) as given in
Table~\ref{testsFunction}. Error in the numerical solution is computed in the
same manner as for meshfree solutions, as per Equation~\eqref{errirInf}.
The same point clouds that will later be used for
assessing point cloud quality metrics are used here as well.
Starting with a uniform $40 \times 40$ Cartesian point cloud with $h=0.025$ and $n=1600$, the points are
perturbed to obtain $1001$ different point clouds, as explained in
Section~\ref{nodSetGeneration}. Each point cloud is tessellated
to generate a corresponding mesh using a Delaunay triangulation, on which we evaluate mesh quality using the inverse aspect ratio.
On each mesh, we compute a numerical solution using the Finite Element Method with linear
shape functions, with the
MATLAB PDE Toolbox \cite{MATLAB_PDE_Toolbox}.

\begin{figure}[!htp]
	\centering
	\includegraphics[width=0.6\textwidth]{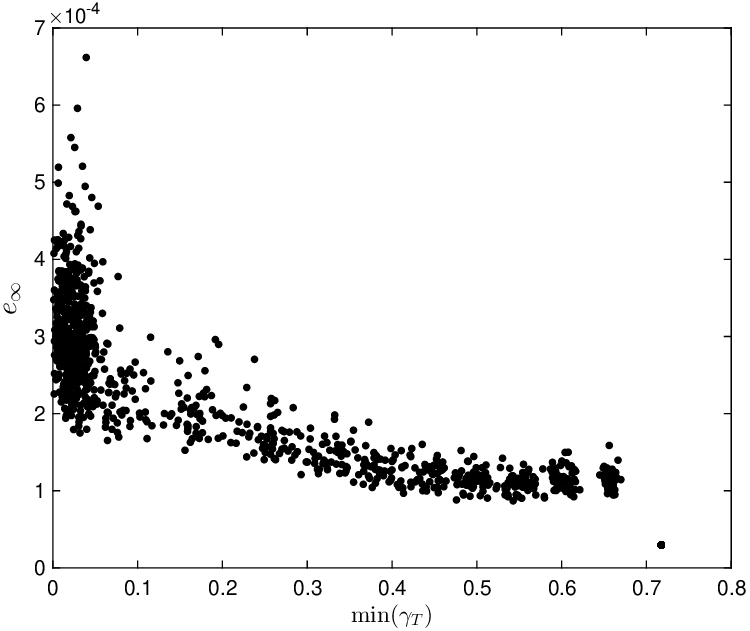} 
\caption{Validation of methodology. Poisson problem in $\mathbb{R}^2$ solved by the Finite Element Method: Variation of the
error \(e_{\infty}\) versus the mesh-quality metric given by the global inverse aspect ratio  $\text{min}(\gamma_{T})$. Each dot represents
the numerical result on a different mesh, each of which is a Delaunay triangulation of the same number of points.
Spearman rank correlation $\rho = 0.872$.}
	\label{aspectRatioFEM}
\end{figure}

Numerical results for all $1001$ meshes are collected in Figure~\ref{aspectRatioFEM} that
shows the variation of the numerical error $e_{\infty}$ with respect to the global inverse aspect ratio metric $\text{min}(\gamma_{T})$. 
The figure reveals a strong counter-monotonic relationship, where higher errors
correspond to lower inverse aspect ratios. This is expected since the aspect ratio, in both its standard and inverse (reciprocal) forms, 
is widely and successfully used as a mesh quality metric. However, the linear correlation between these quantities
is obviously rather weak.


%

While a visual inspection of the error-metric plot helps to identify general trends and directional
relationships, it lacks the precision needed for quantitative assessment. Given the large number
of metrics considered in this work, visual analysis alone is impractical. To systematically evaluate
how quality metrics correlate with the numerical error, we therefore use monotonicity measures (see
Section~\ref{sec:monotonicity}). As a baseline for comparison, we note that in the above finite element test case
the Spearman rank correlation yields $\rho = 0.872$,
and the Kendall's Tau measure gives $\tau = 0.686$. These values provide a reference to assess
the performance of meshfree quality metrics.

\subsubsection*{Quality of meshfree point clouds}

In order to complete the validation of our methodology and justify the selection of a
monotonicity measure for our main tests, we generate the same type of plots as the one in Figure~\ref{aspectRatioFEM} for the error of
the meshfree solution of the Poisson problem versus six global quality metrics with distinctly different behavior, and
verify how well the five monotonicity measures of Section~\ref{sec:monotonicity} predict their performance.


For the same Poisson equation in $[0,1]^2$ with analytical solution $u_1$ and the same point clouds with $h=0.025$
and $n=1600$ we
now compute the error $e_\infty$ of the discrete solution obtained by the meshfree generalized
finite difference method specified in Section~\ref{numericalMethods}, with the polynomial order $p=4$ and the supports
\(\mathbf{S}_{i}\) consisting of  all points within the distance $\lambda = 0.0575$ of $\mathbf{x}_i$.
This corresponds to
$\varepsilon= 2.3$ in the parameter study of Section~\ref{sec:POrder}.

\begin{figure}[!htp]
	\centering
		\begin{subfigure}[!htp]{0.45\textwidth}
		\centering
		\includegraphics[width=\textwidth]{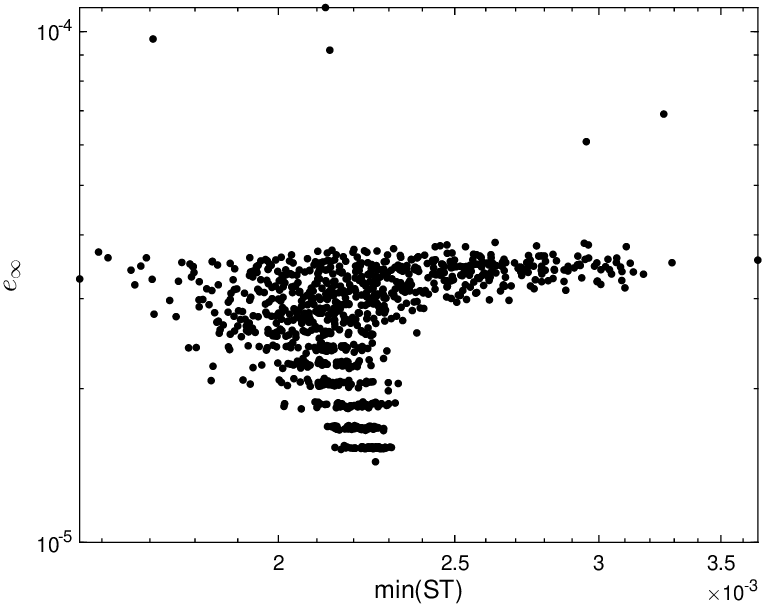}
		\caption{Error versus $\text{min(ST)}$}
		\label{minST}
	\end{subfigure}
	\begin{subfigure}[!htp]{0.45\textwidth}
		\centering
		\includegraphics[width=\textwidth]{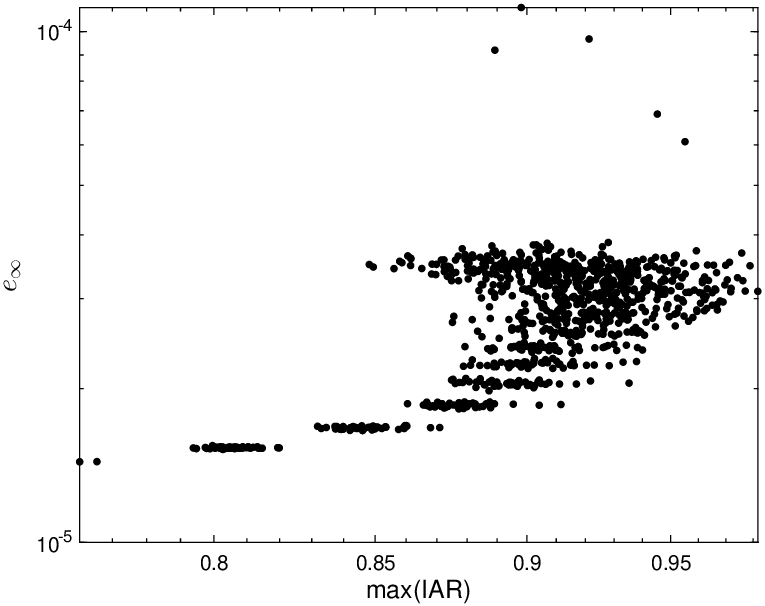}
		\caption{Error versus $\text{max(IAR)}$}
		\label{maxIAR}
	\end{subfigure}
	\\
	\begin{subfigure}[!htp]{0.45\textwidth}
		\centering
		\includegraphics[width=\textwidth]{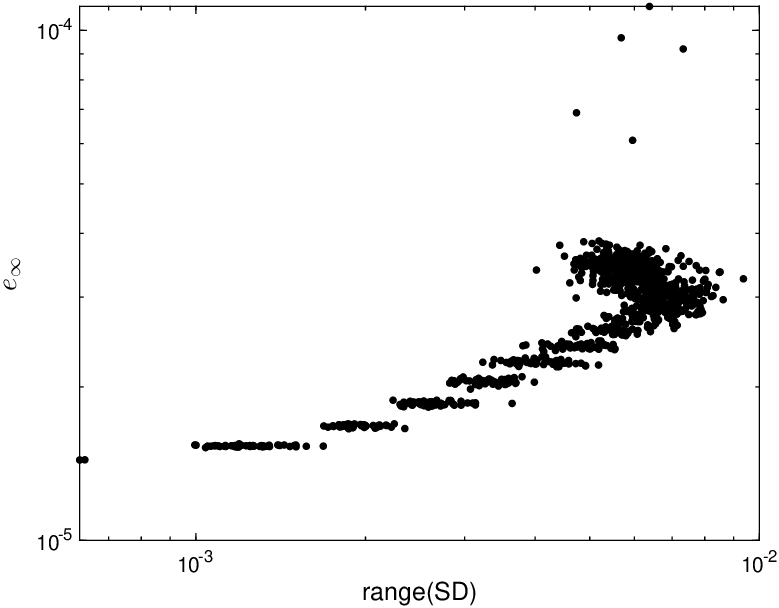}
		\caption{Error versus $\text{range(SD)}$}
		\label{rangeSD}
	\end{subfigure}
	\begin{subfigure}[!htp]{0.45\textwidth}
		\centering
		\includegraphics[width=\textwidth]{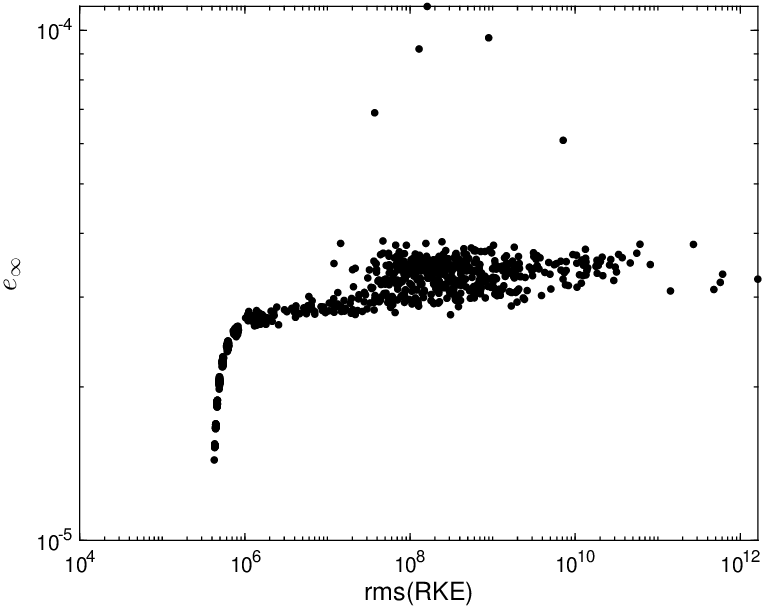}
	\caption{Error versus $\text{rms(RKE)}$}
		\label{rmsRKE}
	\end{subfigure}
	\\
	\begin{subfigure}[!htp]{0.45\textwidth}
		\centering
		\includegraphics[width=\textwidth]{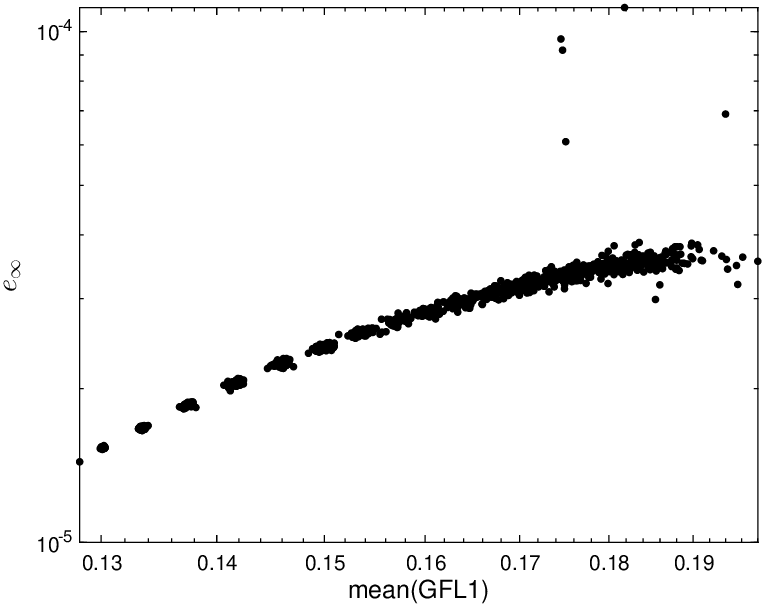}
		\caption{Error versus $\text{mean(GFL1)}$}
		\label{meanGFL1}
	\end{subfigure}
	\begin{subfigure}[!htp]{0.45\textwidth}
		\centering
		\includegraphics[width=\textwidth]{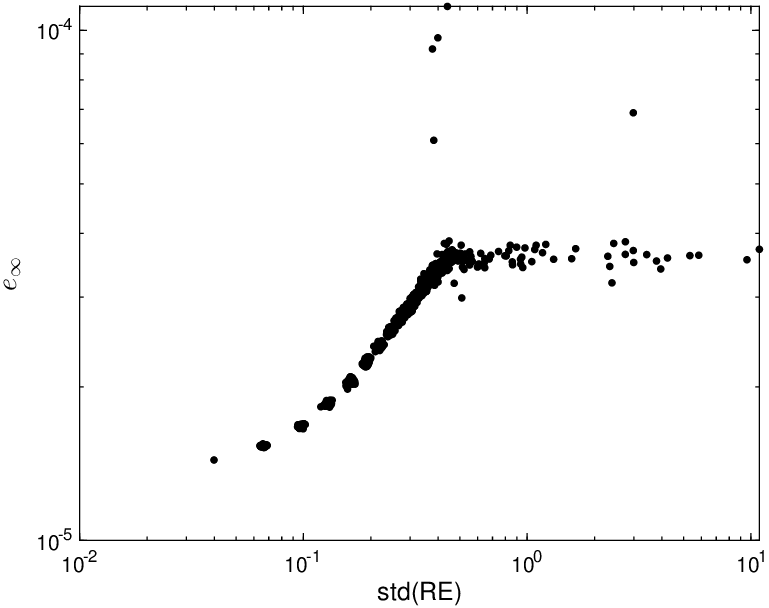}
		\caption{Error versus $\text{std(RE)}$}
		\label{stdRE}
	\end{subfigure}
\caption{Validation of methodology. Poisson problem in $\mathbb{R}^2$ solved by a meshfree method:
Variation of the error \(e_{\infty}\) versus six global quality metrics.
In each figure, every dot represents the numerical result on a different point cloud. The top two figures show a poor
correlation between error and metric values. In contrast, the two middle figures reveal a moderate correlation, whereas the bottom two figures demonstrate a strong monotonic dependence. This
correlation is quantified in Table~\ref{tablemonotoncMeasurevalidation}.}
	\label{monotoncMeasurevalidation}
\end{figure}

Figure~\ref{monotoncMeasurevalidation} shows the variation of the numerical error with respect to six
selected quality metrics from our pool described in Table~\ref{tab:qualityMetrics} and
Appendix~\ref{qualityMetricDetails}.  The first two metrics, $\text{min(ST)}$ (Figure~\ref{minST}) 
and $\text{max(IAR)}$ (Figure~\ref{maxIAR}), exhibit a poor correlation with the error. 
The next one, $\text{range(SD)}$ (Figure~\ref{rangeSD}), 
shows a moderate correlation, while the performance of $\text{rms(RKE)}$ (Figure~\ref{rmsRKE}) is similar to that of the 
inverse aspect ratio in the mesh test of Figure~\ref{aspectRatioFEM}, and the last two metrics, 
$\text{mean(GFL1)}$ and $\text{std(RE)}$, demonstrate a very strong monotonic correlation with the error.

Table~\ref{tablemonotoncMeasurevalidation} provides the absolute values of all five monotonicity measures
for each of the quality metrics used in Figure~\ref{monotoncMeasurevalidation}. The results confirm that Spearman’s rank
correlation coefficient $\rho$ and Kendall's Tau coefficient $\tau$ 
provide the most reliable quantification of the monotonic relationship between quality metrics and numerical error, 
in agreement with the observed trends in Figure~\ref{monotoncMeasurevalidation}. 
Both measures consistently identify $\text{mean(GFL1)}$ and $\text{std(RE)}$ as strong indicators of numerical error, whereas $\text{min(ST)}$ and $\text{max(IAR)}$ are not. 
Furthermore, both measures classify $\text{range(SD)}$ as a moderate indicator, and give close numerical values for 
$\text{rms(RKE)}$ and inverse aspect ratio in the mesh test, again consistent with the figures.
%
In contrast to this, the other three monotonicity measures fail to capture this relationship accurately. The Pearson
correlation coefficient $r$, which is designed to detect linear relationships, does not reflect the high
monotonic correlation seen in Figures~\ref{rmsRKE} and \ref{stdRE}.
Similarly, both the Minimum Covariance Determinant (MCD), and the Monotonicity Coefficient (MC) provide unreliable quantifications. 
The MCD measure incorrectly suggests that ST has a stronger correlation with the error than SD, while the MC measure 
wrongly identifies IAR as a stronger monotonic metric than RE, and sees RKE as the weakest metric,
clearly contradicting empirical observations.   

\begin{table}[!htp]
	\centering
	\caption{Validation of methodology. Performance of monotonicity measures to quantify the relationship
		between	point cloud quality metrics and numerical error for the data presented in Figure~\ref{monotoncMeasurevalidation}.
	Monotonicity measures:  Spearman's rank correlation $\rho$, Kendall's Tau $\tau$, Pearson's correlation
	coefficient $r$, Minimum Covariance Determinant MCD, and Monotonicity Coefficient MC. The table shows that
	$\rho$ and $\tau$ capture the monotonic correlation between the metrics and the error well, while the other three do not.}
	\label{tablemonotoncMeasurevalidation}
	\begin{tabular}{lccccccc}
		\toprule
		\textbf{Metric } & \textbf{Map} & \textbf{Type} & $\rho$ & $\tau$ & $r$ & MCD & MC   \\
		\midrule
		ST & min & system & 0.204 & 0.098 & 0.222&0.491 &0.274 \\
		IAR & max & mesh & 0.504 & 0.378 & 0.601&0.106&0.728  \\
		SD & range & geometry & 0.628 & 0.460 &0.728&0.329 & 0.831 \\	
		RKE & rms & energy	 & 0.895 & 0.758 & 0.055&0.537 &0.228  \\
		GFL1 & mean & system & 0.988 & 0.922 & 0.875&0.999&0.984  \\
		RE & std & system & 0.990 & 0.926 & 0.329&0.998&0.617  \\
		\bottomrule
	\end{tabular}
\end{table}

This trend, where \(\rho\) and \(\tau\) consistently capture the monotonic relationship between error and metric values
well, while the other three monotonicity measures fail, holds across all global metrics considered in Section~\ref{metrics}. For
brevity, we reported results only for the above selection of six global metrics.
Based on this, we conclude that only \(\rho\)  and \(\tau\) are suitable as monotonicity measures in our setting,
and that they are effective in identifying quality metrics that reliably correlate with numerical error.
Throughout our experiments, we observed that \(\rho\)  and \(\tau\) gave very similar results and perform equally well
for our purpose,
henceforth we only report the values of \(\rho\) in the coming sections.



\section{Evaluation of quality metrics}\label{eval}

We now present numerical results to determine which metrics serve as good indicators for point cloud quality.
Section~\ref{parameterStudy} evaluates their
performance for the Poisson problem in $\mathbb{R}^2$ under different
parameters of the meshfree numerical method by varying the point cloud size, support size, and polynomial order. This is
followed by the $\mathbb{R}^3$ Poisson test case in Section~\ref{sec:3D}, and
the hyperbolic problem in Section~\ref{sec:HyperbolicResults}. The results are summarized in Section~\ref{sec:Discussion}.

\subsection{Parameter sensitivity of quality metrics}\label{parameterStudy}


A robust and reliable quality metric should consistently
perform well across different numerical settings. To assess this, we conduct a parameter study by
varying three fundamental parameters of a meshfree discretization:
\begin{itemize}
\setlength{\itemsep}{1pt}
\setlength{\parskip}{0pt}
\setlength{\parsep}{0pt}
    \item \emph{Resolution}: the number of points in the domain, controlling the point cloud size,
	see Section~\ref{sec:Resolution}.
    \item \emph{Support size}: the size of the support \(\mathbf{S}_i\) used for local computations,
    see Section~\ref{sec:Support}.
    \item \emph{Polynomial order}: governing the degree of polynomials used in the numerical approximation,
    see Section~\ref{sec:POrder}.
\end{itemize}
We analyze the performance of all quality metrics introduced in Section~\ref{metrics} across these scenarios
while maintaining the same PDE setting as in Section~\ref{sec:Validation_Results}, the Poisson problem on a unit
square in \(\mathbb{R}^2\), with both analytical solutions  \( u_1 \) and \( u_2 \) of Table~\ref{testsFunction}.

Henceforth, for each local metric listed in Table~\ref{tab:qualityMetrics}, we only report the global metric
that exhibits the highest correlation with the numerical error, while all aggregation maps have been evaluated for each local
metric. 

\subsubsection{Parameter study: Resolution}
\label{sec:Resolution}

We start the parameter study by varying the point cloud resolution, choosing the number of points
 \(n =\) 100, 400, 1600, 6400 and 25600. For each of the five resolutions,
we perturb the Cartesian grid to generate a total of $1001$ point clouds in $[0,1]^2$, as explained in
Section~\ref{nodSetGeneration}.

For a given point cloud $\mathbf{X}$, at each point \(\mathbf{x}_i \in \mathbf{X}\), the support
\(\mathbf{S}_{i}\) consists of all points within a distance $\lambda$ of
$\mathbf{x}_i$, given by
%
\begin{equation}
	\lambda = \varepsilon \, n^{-1/2},
	\label{range}
\end{equation}
where \(\varepsilon\) is a scaling factor. 
In meshfree literature, the support domain is defined either by a fixed number of nearest neighbors
or by selecting all neighbors within a prescribed distance. We adopt the latter approach here, in order to
remain consistent with the MESHFREE code used later for hyperbolic problems (Section~\ref{sec:HyperbolicResults}).
Note that MESHFREE employs this strategy because of its advantages in simplifying domain decomposition in parallel simulations.
Empirically, we also observed that choosing a fixed number of neighbors instead does not
significantly affect the results.

In the resolution tests we take \(\varepsilon = 2.1\), which
results in approximately 25 neighbors for interior points.
The polynomial order is set to $p = 3$, that is, the Laplacians of all
monomials of degree two and less are differentiated exactly in Equation~\eqref{diffFormula}.

For each local quality metric in Table~\ref{tab:qualityMetrics}, we obtain seven global metrics by considering
all the aggregation maps mentioned in Section~\ref{locaAndGlobalDeff}. However, as mentioned above, for the sake of
brevity, we only report the results of the global metric that exhibits the strongest correlation with
numerical error.
To identify the best global metrics, we select, for each local metric, the one with the highest average
Spearman's rank correlation \(\bar{\rho}\) across all five point cloud resolutions.
Tables~\ref{tab:spearma-setsize-u1} and~\ref{tab:spearma-setsize-u2} present the selected global quality
metrics for the \(\mathbb{R}^2\) Poisson problem with analytical solutions \(u_{1}\) and \(u_{2}\) of
Table~\ref{testsFunction},
respectively. In both tables, the metrics are sorted in descending order of \(\bar{\rho}\).

\begin{table}[!ht]
	\centering
	\caption{Parameter study: varying resolution for 2D Poisson problem with analytical solution \( u_1 \).
	Spearman's rank correlation coefficient \(\rho\) is given for the resolutions determined by the number
	\(n\) of points in the point clouds.
 The columns denoted \(\bar{\rho}\) and $\bar\sigma$ give the
	averages of $\rho$ and of the performance score $\sigma$  across all five resolutions.
	The polynomial order is fixed at \( p = 3 \) and the support scaling factor at \( \varepsilon = 2.1 \). }
	\label{tab:spearma-setsize-u1}
		\begin{tabular}{lccccccccc}
		\toprule
					\textbf{Metric } & \textbf{Map} & \textbf{Type} &\multicolumn{5}{c}{\(\rho\)}& \textbf{\(
					\bar{\rho}\)} & $\bar\sigma$\\
		\cmidrule(lr){4-8}
		& & & \textbf{$
			n =100$} & \textbf{$n =400$} & \textbf{$
			n =1600$} & \textbf{$n =6400$} & \textbf{$
			n =25600$}& \\
		\midrule
	GFL2 & rms & system      & 0.774 & 0.820 & 0.836 & 0.810 & 0.819 & 0.812 & 1 \\
	GFL1 & std & system      & 0.777 & 0.819 & 0.834 & 0.808 & 0.819 & 0.812 & 1 \\
	RE  & std & system      & 0.770 & 0.820 & 0.834 & 0.809 & 0.820 & 0.811 & 2 \\
	DJM  & std & mesh        & 0.773 & 0.817 & 0.835 & 0.809 & 0.819 & 0.811 & 2 \\
	IMPF & mean & equilibrium& 0.771 & 0.819 & 0.833 & 0.809 & 0.820 & 0.810 & 3 \\
	DW   & mean & system     & 0.771 & 0.816 & 0.834 & 0.811 & 0.819 & 0.810 & 3 \\
	FD   & std  & geometry       & 0.764 & 0.817 & 0.836 & 0.809 & 0.819 & 0.809 & 4 \\
	IME  & median & energy   & 0.758 & 0.811 & 0.834 & 0.807 & 0.820 & 0.806 & 7 \\
	RKPF & median & equilibrium & 0.752 & 0.816 & 0.834 & 0.807 & 0.821 & 0.806 & 7 \\
	RKE & median & energy   & 0.755 & 0.812 & 0.833 & 0.808 & 0.820 & 0.806 & 7 \\
	RLR  & mean & geometry   & 0.752 & 0.814 & 0.833 & 0.802 & 0.820 & 0.804 & 9 \\
	MCD  & std & mesh        & 0.742 & 0.811 & 0.833 & 0.807 & 0.819 & 0.802 & 11 \\
	IAR   & std & mesh        & 0.749 & 0.808 & 0.815 & 0.804 & 0.812 & 0.798 & 15 \\
	LKPF  & mean & equilibrium& 0.708 & 0.810 & 0.832 & 0.810 & 0.819 & 0.796 & 17 \\
	SD   & rms & geometry    & 0.737 & 0.804 & 0.819 & 0.795 & 0.814 & 0.794 & 19 \\
	PD   & mean & geometry   & 0.737 & 0.805 & 0.818 & 0.794 & 0.814 & 0.794 & 19 \\
	CN   & median & system   & 0.716 & 0.806 & 0.823 & 0.804 & 0.818 & 0.793 & 20 \\
	WKPF & mean & equilibrium& 0.693 & 0.808 & 0.832 & 0.809 & 0.819 & 0.792 & 21 \\
	UR   & median & geometry     & 0.728 & 0.806 & 0.817 & 0.796 & 0.814 & 0.792 & 21 \\
	CPF  & mean & equilibrium& 0.679 & 0.795 & 0.827 & 0.803 & 0.808 & 0.782 & 31 \\
	ALR  & std & geometry    & 0.625 & 0.797 & 0.829 & 0.810 & 0.818 & 0.776 & 37 \\
	ST   & median & system   & 0.611 & 0.780 & 0.828 & 0.806 & 0.820 & 0.769 & 44 \\
	WKE  & std & energy      & 0.556 & 0.799 & 0.830 & 0.808 & 0.818 & 0.762 & 51 \\
	CE   & std & energy      & 0.537 & 0.772 & 0.828 & 0.809 & 0.818 & 0.753 & 60 \\
	LKE  & std & energy      & 0.460 & 0.782 & 0.824 & 0.808 & 0.818 & 0.739 & 74 \\
		\bottomrule
	\end{tabular}
\end{table}
\begin{table}[!ht]
	\centering
	\caption{Parameter study: varying resolution for 2D Poisson problem with analytical solution \( u_2 \).
	Spearman's rank correlation coefficient \(\rho\) is given for the resolutions determined by the number
	\(n\) of points in the point clouds.
 The  columns denoted \(\bar{\rho}\) and $\bar\sigma$ give the
	averages of $\rho$ and of the performance score $\sigma$ across all five resolutions.
	The polynomial order is fixed at \( p = 3 \) and the support scaling factor at \( \varepsilon = 2.1 \).}
	\label{tab:spearma-setsize-u2}
		\begin{tabular}{lccccccccc}
			\toprule
			\textbf{Metric } & \textbf{Map} & \textbf{Type} &\multicolumn{5}{c}{\(\rho\)}& \textbf{\( \bar{\rho}\)} & $\bar\sigma$\\
			\cmidrule(lr){4-8}
			& & & \textbf{$
				n =100$} & \textbf{$n =400$} & \textbf{$
				n =1600$} & \textbf{$n =6400$} & \textbf{$
				n =25600$}& \\
			\midrule
	GFL2 & rms    & system      & 0.805 & 0.870 & 0.857 & 0.837 & 0.793 & 0.832 & 1 \\
	GFL1 & std    & system      & 0.798 & 0.873 & 0.860 & 0.838 & 0.791 & 0.832 & 1 \\
	RE  & std    & system      & 0.796 & 0.872 & 0.859 & 0.838 & 0.792 & 0.831 & 2 \\
	FD   & rms    & geometry        & 0.798 & 0.872 & 0.857 & 0.837 & 0.792 & 0.831 & 2 \\
	DJM  & std    & mesh        & 0.793 & 0.874 & 0.858 & 0.837 & 0.791 & 0.831 & 2 \\
	IMPF & mean   & equilibrium & 0.786 & 0.868 & 0.860 & 0.837 & 0.792 & 0.829 & 4 \\
	DW   & mean   & system      & 0.783 & 0.870 & 0.857 & 0.837 & 0.792 & 0.828 & 6 \\
	IME  & median & energy      & 0.780 & 0.865 & 0.854 & 0.835 & 0.792 & 0.825 & 8 \\
	RKE & median & energy      & 0.776 & 0.865 & 0.854 & 0.836 & 0.793 & 0.825 & 8 \\
	RKPF & median & equilibrium & 0.777 & 0.864 & 0.854 & 0.836 & 0.793 & 0.825 & 8 \\
	RLR  & mean   & geometry    & 0.780 & 0.862 & 0.851 & 0.836 & 0.792 & 0.824 & 9 \\
	MCD  & std    & mesh        & 0.770 & 0.866 & 0.854 & 0.836 & 0.793 & 0.824 & 9 \\
	IAR   & std    & mesh        & 0.776 & 0.860 & 0.840 & 0.816 & 0.777 & 0.814 & 19 \\
	CN   & median & system      & 0.748 & 0.841 & 0.845 & 0.831 & 0.794 & 0.812 & 21 \\
	SD   & rms    & geometry    & 0.763 & 0.845 & 0.840 & 0.823 & 0.788 & 0.812 & 21 \\
	PD   & mean   & geometry    & 0.764 & 0.845 & 0.839 & 0.822 & 0.788 & 0.812 & 21 \\
	WKPF & median & equilibrium & 0.721 & 0.853 & 0.855 & 0.837 & 0.792 & 0.811 & 22 \\
	LKPF  & mean   & equilibrium & 0.718 & 0.854 & 0.855 & 0.836 & 0.792 & 0.811 & 22 \\
	UR   & median & geometry        & 0.760 & 0.839 & 0.833 & 0.819 & 0.788 & 0.808 & 25 \\
	CPF  & mean   & equilibrium & 0.724 & 0.846 & 0.844 & 0.829 & 0.785 & 0.806 & 27 \\
	WKE  & std    & energy      & 0.641 & 0.853 & 0.853 & 0.837 & 0.792 & 0.795 & 38 \\
	ST   & median & system      & 0.671 & 0.826 & 0.849 & 0.837 & 0.792 & 0.795 & 38 \\
	ALR  & std    & geometry    & 0.634 & 0.849 & 0.856 & 0.836 & 0.792 & 0.794 & 39 \\
	LKE  & std    & energy      & 0.576 & 0.840 & 0.852 & 0.836 & 0.792 & 0.779 & 54 \\
	CE   & std    & energy      & 0.563 & 0.830 & 0.854 & 0.836 & 0.792 & 0.775 & 58 \\
		\bottomrule
	\end{tabular}
\end{table}

We observe that in both tables \texttt{rms}, \texttt{std} and \texttt{mean} feature as the best aggregation maps for better
performing metrics, while \texttt{median} occasionally appears as the best aggregation map for some metrics in the lower parts of the tables.
However, the best map varies depending on the local metric. 

An important observation is that the best monotonicity correlation $\rho$ achieved over all metrics significantly varies with
the parameter $n$ and between $u_1$ and  $u_2$. In particular, for $n=400$ and $u_2$ we get the highest value
$\rho=0.870$, whereas  $n=100$ and $u_1$ allows at most $\rho=0.774$. However, a difference of about 0.1 between the values of
$\rho$ may already  mean a big performance difference, compare for example Figures~\ref{rmsRKE} and \ref{meanGFL1}, where the
difference in $\rho$ is just 0.092. A consequence of this fact is that it is impossible to provide a single tolerance value
for $\rho$ that would allow to distinguish between good and bad performance of quality metrics. The reason for this behavior
seems to be that the influence of the point cloud quality on the numerical error of the solution is more pronounced in some
settings than in others.

Nevertheless, the data in Tables~\ref{tab:spearma-setsize-u1} and~\ref{tab:spearma-setsize-u2} show that we may always compare
the values of $\rho$ for different metrics applied to any fixed set of point clouds with a particular setting of the
parameters of the meshfree method. In order to demonstrate this, we introduce the \emph{performance score} $\sigma\ge 1$ of a
quality metric defined as
\begin{equation}\label{sigma}
\sigma=10^3(\rho_{\max}-\rho)+1,
\end{equation}
where $\rho$ is Spearman's correlation coefficient of the current metric, and $\rho_{\max}$ is the maximum of $\rho$ over all
metrics tested in a particular experiment where all parameters of the meshfree numerical method are fixed. Thus,
the coefficient $\rho$ of the metric in question is by
$10^{-3}(\sigma-1)$ lower than the best $\rho=\rho_{\max}$ achieved by any metric in the experiment. 
The lower $\sigma$ is, the better the performance of the metric.


The main purpose of $\sigma$ is to allow
a quantitative comparison of the performance of the metrics across the experiments with different parameters. The last columns
of  Tables~\ref{tab:spearma-setsize-u1} and~\ref{tab:spearma-setsize-u2} give the average $\bar\sigma$ of the
performance scores  $\sigma$ in the experiments presented in the table, rounded to the next integer. Ignoring possible  small
inaccuracies due to inconsistent rounding, we compute $\bar\sigma$ as $10^3(\bar\rho_{\max}-\bar\rho)+1$, 
where $\bar\rho$ is the average Spearman's rank correlation rounded to the third digit after the comma, as presented in the 
penultimate columns of the tables, and $\bar\rho_{\max}$ is its maximum achieved for the best performing metric, which therefore
gets the highest possible score $\bar\sigma=1$.


The usefulness of the performance score stems from its remarkable stability, at least for better performing metrics,
over all our parameter studies. In particular, the difference of the values of $\bar\sigma$ between Tables~\ref{tab:spearma-setsize-u1}
and~\ref{tab:spearma-setsize-u2} does not exceed 6 for all metrics but  WKE and LKE whose performance is
weak anyway. This means that the average monotonicity measure \(\bar{\rho}\) remains consistent across the two analytical
solutions $u_1$ and $u_2$, with the variation not exceeding $10^{-3}\bar\sigma=0.006$ after the above adjustment by
$\bar\rho_{\max}$. Hence, the best-performing metrics for
one analytical solution are also the best for the other. This suggests that the performance is
largely independent of the solution, indicating robustness and reliability of the suggested assessment approach.
Therefore, we consider in this section only one analytical solution for each of the numerical tests to follow, with
results for the second solution provided in Appendix~\ref{sec:AppResults}.

The best performing quality metrics, listed in the upper rows of Tables~\ref{tab:spearma-setsize-u1} and
\ref{tab:spearma-setsize-u2}, also exhibit consistency across different resolutions. In particular, for each of
the first 10 metrics (which are the same in both tables)  the difference between the best and the worst scores $\sigma$ over
all resolutions in both Tables~\ref{tab:spearma-setsize-u1} and
\ref{tab:spearma-setsize-u2} do not exceed 29, or even 10 if the standing somewhat apart case $n=100$ is excluded.
This means that the adjusted differences in $\rho$ are at most $0.029$ or $0.01$, respectively.
This further supports the robustness of many metrics in the upper parts of the tables.
A more detailed discussion of the results is presented in Section~\ref{sec:Discussion}.

Based on these observations, we will evaluate the performance of the metrics by comparing their scores
$\sigma$ in our tests, with a smaller score indicating a better performance.

\subsubsection{Parameter study: Support size}
\label{sec:Support}

We now examine the performance of the quality metrics under varying support sizes. For a given point
cloud \(\mathbf{X}\), the support \(\mathbf{S}_{i}\) of each point \(\mathbf{x}_i \in \mathbf{X}\) is
determined by proximity, following Equation~\eqref{range}. As the coefficient \(\varepsilon\)
in \eqref{range} increases, the support size increases.  Here, we consider four values of
\(\varepsilon\): 1.8, 2.1, 2.5, and 2.8 which  correspond approximately to 20, 25, 35, and 45 points  in the
support of interior points. Naturally, points near the boundary will have fewer neighboring points.

We proceed in a similar way to the resolution variation study above, reporting only the
best performing global metric for each local quality metric. The selected aggregation map is determined based on the
highest average Spearman's  rank correlation across all four support sizes considered.
Table~\ref{tab:spearma-subsetsize-u1} presents the monotonic behavior
of the global metrics, for the same analytical solution \( u_1 \) as above. The metrics in the
table are listed in descending order of their average monotonicity correlation.
The corresponding results for the analytical function \( u_2 \) are listed in Table~\ref{tab:spearma-subsetsize-u2}
in Appendix~\ref{sec:AppResults}.

The results show that the monotonic dependence trends for global metrics remain consistent with those
observed in the point cloud resolution study above, reinforcing the robustness and reliability of the best
performing quality metrics across different scenarios. Other trends observed in Section~\ref{sec:Resolution}
also hold here. Moreover, the best and worst performing metrics remain largely the same,
with only slight variations in ordering. For some local metrics, such as GFL1, the best aggregation map for obtaining
the global metric differs from that in the resolution study, though the difference in performance
between both aggregation maps (\texttt{std}, \texttt{rms}) is minimal. This further
suggests that, for many local metrics, multiple aggregation maps  are  viable candidates for computing global metrics,
see Table~\ref{Hperform} and a discussion in Section~\ref{sec:Discussion}.

%
\begin{table}[!ht]
	\centering
	\caption{Parameter study: varying support size with analytical solution \( u_1 \) (see
		Table~\ref{testsFunction}).
	Spearman's rank correlation values \(\rho\) between each quality metric and numerical error for
	different support sizes. The support size is given by $\varepsilon$, which represents
	the scaling factor (see Equation~\eqref{range}).  The  columns denoted \(\bar{\rho}\) and $\bar\sigma$ give the
	averages of $\rho$ and of the performance score $\sigma$ across all four scaling factors. The other parameters are set to
	polynomial order \( p = 3 \) and $n = 1600$ points in the domain.} 
	\label{tab:spearma-subsetsize-u1}
		\begin{tabular}{lcccccccc}
			\toprule
						\textbf{Metric} & \textbf{Map} & \textbf{Type} & \multicolumn{4}{c}{\(\rho\)}& \textbf{\(
						\bar{\rho}\)} & $\bar\sigma$\\
			\cmidrule(lr){4-7}
			& & & \textbf{$\varepsilon =1.8$}&\textbf{$\varepsilon =2.1$}&\textbf{$\varepsilon = 2.5$}
			&\textbf{$\varepsilon=2.8$} & \\
			\midrule
		FD   & rms    & geometry        & 0.839 & 0.837 & 0.874 & 0.885 & 0.859 & 1 \\
		GFL1 & rms    & system      & 0.839 & 0.836 & 0.873 & 0.886 & 0.859 & 1 \\
		RE  & rms    & system      & 0.838 & 0.836 & 0.874 & 0.885 & 0.858 & 2 \\
		GFL2 & rms    & system      & 0.838 & 0.836 & 0.872 & 0.884 & 0.857 & 3 \\
		DJM  & std    & mesh        & 0.836 & 0.835 & 0.872 & 0.882 & 0.856 & 4 \\
		IME  & median & energy      & 0.834 & 0.834 & 0.872 & 0.884 & 0.856 & 4 \\
		RKE & median & energy      & 0.834 & 0.833 & 0.871 & 0.883 & 0.856 & 4 \\
		RKPF & median & equilibrium & 0.834 & 0.834 & 0.871 & 0.883 & 0.855 & 5 \\
		LKPF  & rms    & equilibrium & 0.834 & 0.832 & 0.874 & 0.880 & 0.855 & 5 \\
		DW   & mean   & system      & 0.836 & 0.834 & 0.870 & 0.881 & 0.855 & 5 \\
		MCD  & std    & mesh        & 0.834 & 0.833 & 0.871 & 0.882 & 0.855 & 5 \\
		RLR  & mean   & geometry    & 0.831 & 0.833 & 0.871 & 0.884 & 0.855 & 5 \\
		IMPF & median & equilibrium & 0.833 & 0.834 & 0.870 & 0.881 & 0.855 & 5 \\
		WKPF & mean   & equilibrium & 0.833 & 0.832 & 0.873 & 0.879 & 0.854 & 6 \\
		CPF  & mean   & equilibrium & 0.834 & 0.827 & 0.875 & 0.880 & 0.854 & 6 \\
		WKE  & std    & energy      & 0.832 & 0.830 & 0.866 & 0.870 & 0.849 & 11 \\
		CE   & std    & energy      & 0.832 & 0.828 & 0.864 & 0.869 & 0.848 & 12 \\
		ALR  & std    & geometry    & 0.832 & 0.829 & 0.864 & 0.867 & 0.848 & 12 \\
		LKE  & std    & energy      & 0.831 & 0.824 & 0.862 & 0.867 & 0.846 & 14 \\
		IAR   & std    & mesh        & 0.824 & 0.815 & 0.858 & 0.870 & 0.842 & 18 \\
		CN   & median & system      & 0.809 & 0.823 & 0.862 & 0.872 & 0.842 & 18 \\
		SD   & rms    & geometry    & 0.817 & 0.819 & 0.857 & 0.870 & 0.841 & 19 \\
		PD   & mean   & geometry    & 0.817 & 0.818 & 0.857 & 0.869 & 0.840 & 20 \\
		UR   & median & geometry        & 0.818 & 0.817 & 0.855 & 0.869 & 0.840 & 20 \\
		ST   & max    & system      & 0.802 & 0.785 & 0.834 & 0.824 & 0.811 & 49 \\
		 \bottomrule
	\end{tabular}
\end{table}

Note that the performance of some metrics is significantly different here from Section~\ref{sec:Resolution}.
In particular the performance score $\bar\sigma$ of ALR is 12 in Table~\ref{tab:spearma-subsetsize-u1} and 4 in
Table~\ref{tab:spearma-subsetsize-u2} versus 37 in Table~\ref{tab:spearma-setsize-u1} and 39 in
Table~\ref{tab:spearma-setsize-u2}.  Such metrics that perform well in some parameter studies but poorly in
others are potentially unreliable, unless the scope of their applications can be precisely determined.
We will discuss this in more detail in Section~\ref{sec:Discussion}.

\subsubsection{Parameter study: Polynomial order}
\label{sec:POrder}

As the final parameter dependence study, we now consider the performance of the quality metrics under
varying polynomial order, including higher order meshfree discretizations.
High-order methods improve accuracy for smooth solutions, such as $u_1$ and $u_2$,
but the quality of point clouds remains a critical
factor influencing the numerical error. 

For the Poisson equation on a unit square in \(\mathbb{R}^2\), consider the numerical method of
Section~\ref{numericalMethods} with polynomial
orders $p= 3$, $4$, and $5$ while keeping the point cloud resolution fixed at
\(n = 1600\). This corresponds to quadratic, cubic, and quartic exactness,
respectively. The support at each point is still determined according to
Equation~\eqref{range}. For higher-order polynomials (\( p > 3 \)), larger support sizes are required in order to
obtain reliable numerical differentiation formulas \eqref{diffFormula}.
To this end, the coefficient \(\epsilon\) in Equation~\eqref{range} is set to 2.1, 2.3, and 2.8 for polynomial orders 3, 4, and 5,
respectively.

Table~\ref{tab:spearma-porder-u1} presents the monotonic correlation
of global quality metrics for analytical solution \( u_1 \). As in the
previous parameter studies, we only report the best global metric for each local metric, determined by
the highest average Spearman's rank monotonicity measure across all polynomial orders considered.
The table shows that most metrics serve as better indicators of point cloud quality at higher
polynomial orders.
The corresponding results for the analytical function \( u_2 \) are listed in Table~\ref{tab:spearma-porder-u2}
in Appendix~\ref{sec:AppResults}.

\begin{table}[!ht]
	\centering
	\caption{Parameter study: varying polynomial order with analytical solution \( u_1 \) (see
		Table~\ref{testsFunction}).
	Spearman's rank correlation values \(\rho\) between each quality metric and numerical error for
	different polynomial orders. The polynomial order is given by $p$.
	The  columns denoted \(\bar{\rho}\) and $\bar\sigma$ give the
	averages of $\rho$ and of the performance score $\sigma$ across all polynomial orders. 
	The other parameters are set to scaling factors of 2.1, 2.3, and 2.8 for polynomial orders $p = 3, 4,$ and \(5\),
respectively, with $n = 1600$ points in the domain.} 
	\label{tab:spearma-porder-u1}
		\begin{tabular}{lccccccc}
			\toprule
		\textbf{Metric} & \textbf{Map} & \textbf{Type} & \multicolumn{3}{c}{\(\rho\)}&\textbf{\( \bar{\rho}\)}  & $\bar\sigma$\\
		\cmidrule(lr){4-6}
		& & & \textbf{$p =3$} & \textbf{$p=4$} &\textbf{$p =5$} & \\
			\midrule
		FD   & rms    & geometry   & 0.837 & 0.989 & 0.914 & 0.913 & 1 \\
		RE  & std    & system      & 0.834 & 0.990 & 0.914 & 0.913 & 1 \\
		GFL1 & mean   & system      & 0.835 & 0.988 & 0.915 & 0.913 & 1 \\
		GFL2 & mean   & system      & 0.836 & 0.987 & 0.915 & 0.913 & 1 \\
		DJM  & std    & mesh        & 0.835 & 0.989 & 0.913 & 0.912 & 2 \\
		IMPF & rms    & equilibrium & 0.834 & 0.989 & 0.912 & 0.912 & 2 \\
		RKPF & median & equilibrium & 0.834 & 0.988 & 0.912 & 0.911 & 3 \\
		RKE & median & energy      & 0.833 & 0.988 & 0.911 & 0.911 & 3 \\
		IME  & median & energy      & 0.834 & 0.988 & 0.911 & 0.911 & 3 \\
		LKPF  & mean   & equilibrium & 0.833 & 0.986 & 0.911 & 0.910 & 4 \\
		WKPF & mean   & equilibrium & 0.832 & 0.985 & 0.910 & 0.909 & 5 \\
		MCD  & std    & mesh        & 0.833 & 0.985 & 0.908 & 0.909 & 5 \\
		RLR  & median & geometry    & 0.833 & 0.985 & 0.908 & 0.909 & 5 \\
		CPF  & mean   & equilibrium & 0.827 & 0.986 & 0.910 & 0.908 & 6 \\
		ALR  & std    & geometry    & 0.830 & 0.984 & 0.904 & 0.906 & 8 \\
		CE   & std    & energy      & 0.828 & 0.983 & 0.904 & 0.905 & 9 \\
		WKE  & std    & energy      & 0.830 & 0.977 & 0.900 & 0.902 & 12 \\
		LKE  & std    & energy      & 0.824 & 0.972 & 0.896 & 0.897 & 17 \\
		IAR   & std    & mesh        & 0.815 & 0.965 & 0.897 & 0.892 & 22 \\
		SD   & rms    & geometry    & 0.819 & 0.958 & 0.887 & 0.888 & 26 \\
		CN   & median & system      & 0.823 & 0.946 & 0.893 & 0.887 & 27 \\
		PD   & mean   & geometry    & 0.818 & 0.958 & 0.886 & 0.887 & 27 \\
		UR   & median & geometry    & 0.817 & 0.951 & 0.878 & 0.882 & 32 \\
		DW   & std    & system      & 0.808 & 0.948 & 0.885 & 0.881 & 33 \\
		ST   & median & system      & 0.828 & 0.757 & 0.914 & 0.833 & 81 \\
		\bottomrule
	\end{tabular}
\end{table}

Through these parameter studies, we observed that the \texttt{median}, \texttt{mean}, \texttt{standard deviation}, and \texttt{root mean square} aggregation maps are effective choices 
for converting local metrics into global metrics, particularly for high performing metrics. 
For these metrics, the correlation between the resulting global metric and the numerical error is only weakly affected by the choice among these aggregation maps.
The best and worst performing metrics again show a similar pattern to that observed earlier.
For example, \textit{equilibrium} based metrics show a stronger monotonic relationship with error
compared to their corresponding \textit{energy} based metrics, consistent with prior findings.

\subsection{Quality metrics in 3D} \label{sec:3D}
We now extend the analysis to three dimensions to assess whether the trends observed in \(\mathbb{R}^2\)
also hold in \(\mathbb{R}^3\).
The same methodology is applied to the Poisson problem on a unit cube in \(\mathbb{R}^3\), using analytical
solutions $ u_3 $ and $ u_4 $, see Table~\ref{testsFunction}.
We consider two resolutions: a coarse point cloud with \(n = 2744\) and a fine point
cloud with \(n = 39304\). As before, for each resolution, $1001$ point clouds are generated
by random perturbations of a Cartesian grid, see Section~\ref{nodSetGeneration}. The support at each point is given by
proximity by considering all points within a distance \(\lambda\) given by
%
\begin{equation}
	\lambda = \varepsilon \, n^{-1/3},
	\label{3Drange}
\end{equation}
where $n$ is the number of points in the point cloud.
We set the scaling factor to \(\varepsilon = 2.1\) and the polynomial order to $p=3$ for quadratic exactness.

We evaluate the quality metrics in \(\mathbb{R}^3\) in a similar way as done above in \(\mathbb{R}^2\),
by analyzing the monotonic dependency between metric values and numerical error. Table~\ref{tab:3d-u1}
presents the monotonic behavior of the global quality metrics for analytical solution
\( u_3 \), ranked by their average Spearman's correlation across both resolutions
considered. As before, only the best performing global metric is reported for each local metric.
The corresponding results for the analytical function \( u_4 \) are listed in Table~\ref{tab:spearma-3d-u2}
in Appendix~\ref{sec:AppResults}.
\begin{table}[!ht]
	\centering
	\caption{Poisson problem in \(\mathbb{R}^3\) with analytical solution \( u_3 \) (see
		Table~\ref{testsFunction}):
	The Spearman's rank correlation \(\rho\) quantifying the monotonic dependence between quality metrics
	and numerical error. The discretization resolution is given by \(n\), representing the
	number of points in the point clouds considered.   The  columns denoted \(\bar{\rho}\) and $\bar\sigma$ give the
	averages of $\rho$ and of the performance score $\sigma$ across both resolutions.
	Other parameters are set to polynomial order \( p = 3 \) and
	scaling factor $\varepsilon = 2.1$.} 
	\label{tab:3d-u1}
		\begin{tabular}{lcccccc}
			\toprule
		\textbf{Metric} & \textbf{Map} & \textbf{Type} & \multicolumn{2}{c}{\(\rho\)} & \textbf{\(\bar{\rho}\)}  & $\bar\sigma$\\
		\cmidrule(lr){4-5}
		& & & \( n = 2744 \) & \( n = 39304 \) & \\
			\midrule
		RE   & std    & system      & 0.955 & 0.968 & 0.962 & 1 \\
		GFL2  & std    & system      & 0.954 & 0.968 & 0.961 & 2 \\
		FD    & std    & geometry    & 0.953 & 0.969 & 0.961 & 2 \\
		DJM   & std    & mesh        & 0.952 & 0.968 & 0.960 & 3 \\
		IMPF  & rms    & equilibrium & 0.951 & 0.968 & 0.960 & 3 \\
		GFL1  & std    & system      & 0.950 & 0.968 & 0.959 & 4 \\
		DW    & mean   & system      & 0.950 & 0.968 & 0.959 & 4 \\
		RKE  & median & energy      & 0.953 & 0.966 & 0.959 & 4 \\
		RKPF  & median & equilibrium & 0.949 & 0.965 & 0.957 & 6 \\
		CN    & median & system      & 0.946 & 0.967 & 0.956 & 7 \\		
		IAR    & std    & mesh        & 0.943 & 0.967 & 0.955 & 8 \\
		IME   & median & energy      & 0.942 & 0.967 & 0.954 & 9 \\
		CPF   & mean   & equilibrium & 0.947 & 0.960 & 0.953 & 10 \\				
		RLR   & mean   & geometry    & 0.935 & 0.962 & 0.948 & 15 \\
		LKPF   & mean   & equilibrium & 0.926 & 0.968 & 0.947 & 16 \\
		WKPF  & mean   & equilibrium & 0.926 & 0.968 & 0.947 & 16 \\
		WKE   & std    & energy      & 0.933 & 0.956 & 0.944 & 19 \\
		MCD   & std    & mesh        & 0.930 & 0.954 & 0.942 & 21 \\
		LKE   & std    & energy      & 0.911 & 0.960 & 0.936 & 27 \\
		UR    & range  & geometry    & 0.929 & 0.937 & 0.933 & 30 \\
		SD    & max    & geometry    & 0.907 & 0.927 & 0.917 & 46 \\
		PD    & std    & geometry    & 0.906 & 0.921 & 0.914 & 49 \\
		CE    & std    & energy      & 0.860 & 0.967 & 0.914 & 49 \\
		ALR   & std    & geometry    & 0.860 & 0.967 & 0.913 & 50 \\
		ST    & range  & system      & 0.865 & 0.887 & 0.876 & 87 \\
		\bottomrule
	\end{tabular}
\end{table}

The tables confirm a similar behavior to that observed in the \(\mathbb{R}^2\) simulations, with similar
ranking of quality metrics. Metrics such as RE, GFL2, FD, DJM, IMPF, and GFL1 continue to be reliable
indicators of error. Minimal variation in \(\bar{\rho}\) is observed across the two resolutions considered
for the best performing metrics, emphasizing their robustness.
Conversely, weaker metrics like CE and ALR
show greater variability between resolutions, indicating their unreliability.

Furthermore, the tables show that each individual quality metric provides a stronger correlation with error at
the higher resolution, sometimes significantly so. This supports the assumption that the quality metrics
perform better whenever the sources of error other than the point cloud quality diminish, since finer point clouds 
better resolve the features of the solution.


\subsection{Quality metrics for hyperbolic PDEs}
\label{sec:HyperbolicResults}

We now extend our analysis to hyperbolic PDEs by considering the advection equation in \(\mathbb{R}^2\).
While the Poisson problem considered above models diffusion dominated phenomena, the advection equation
describes transport dominated processes.
This allows us to evaluate the robustness of quality metrics across different PDE types.
We follow the same methodology as in the elliptic case, analyzing the correlation between quality metrics
and numerical error.
Details about the analytical solution considered, and spatial and temporal discretization are described
in Sections~\ref{sec:Advection} and \ref{numericalMethods}.

The point clouds considered here contain \(n = 1600\) points and are the same as those
used for the \(\mathbb{R}^2\) Poisson problem in Sections~\ref{sec:Support} and~\ref{sec:POrder}, shifted to the square
$[-0.5,0.5]^2$. For each
point \(\mathbf{x}_i\), the support set \(\mathbf{S}_i\) is determined by proximity, as per Equation~\eqref{range}.
The scaling parameter is set to \(\varepsilon = 3.0\), which results in \(\lambda = 0.075\).
Furthermore, the polynomial order is set to 3, i.e. monomials up to degree two are differentiated exactly.

%
%


Table~\ref{spearmanHyper} presents the monotonic correlation of the global quality metrics
with numerical error. As before, only the best-performing global metric is reported for each
local metric, and the metrics are ordered according to a descending Spearman's rank correlation
coefficient $\rho$.
Table~\ref{spearmanHyper} shows that \textit{system matrix} quality metrics GFL1, GFL2, and RE exhibit strong correlation with
numerical error, consistent with the trends observed for the elliptic PDEs above.
Even if they are  designed to measure the accuracy of the local discretization in space, these results
may indicate that these metrics  are also sensitive to the specific stability issues of
the hyperbolic equations tackled by upwinding and limiters at the expense of lower local accuracy. 
However, quality metric DW,
again in this type, shows a much weaker correlation with the numerical error than that observed in the elliptic PDE results.

%
\begin{table}[h!]
	\centering
		\caption{Advection problem in \(\mathbb{R}^2\): The Spearman's rank correlation \(\rho\) quantifying
		the monotonic dependence between quality metrics and numerical error for the hyperbolic PDE,
		and the performance score $\sigma$. Numerical
		parameters are set to $n = 1600$ points in the domain, polynomial order
		\( p = 3 \) and scaling factor \( \varepsilon = 3.0 \).}
	\label{spearmanHyper}
	\begin{tabular}{lcccc}
	\toprule
		\textbf{Metric} & \textbf{Map} & \textbf{Type} &\textbf{\( \rho\)} & $\sigma$\\
	\midrule
DJM   & std    & mesh        & 0.863 & 1   \\
GFL1  & std    & system      & 0.862 & 2   \\
RE   & std    & system      & 0.862 & 2   \\
GFL2  & std    & system      & 0.862 & 2   \\
FD    & rms    & geometry    & 0.861 & 3   \\
RKE  & median & energy      & 0.860 & 4   \\
IME   & median & energy      & 0.860 & 4   \\
IMPF  & mean   & equilibrium & 0.860 & 4   \\
RKPF  & median & equilibrium & 0.859 & 5   \\
RLR   & median & geometry    & 0.856 & 8   \\
MCD   & std    & mesh        & 0.856 & 8   \\
LKPF  & mean   & equilibrium & 0.855 & 9   \\
WKPF  & median & equilibrium & 0.854 & 10   \\
CPF   & mean   & equilibrium & 0.853 & 11  \\
ALR   & std    & geometry    & 0.848 & 16  \\
CE    & std    & energy      & 0.847 & 17  \\
IAR    & std    & mesh        & 0.841 & 23  \\
WKE   & std    & energy      & 0.838 & 26  \\
CN    & median & system      & 0.827 & 37  \\
LKE   & std    & energy      & 0.824 & 40  \\
PD    & min    & geometry    & 0.822 & 42  \\
SD    & rms    & geometry    & 0.818 & 46  \\
UR    & median & geometry    & 0.812 & 52  \\
ST    & median & system      & 0.736 & 128 \\
DW    & mean   & system      & 0.705 & 159 \\
			\bottomrule
		\end{tabular}
	\end{table}
%


In addition, \textit{equilibrium} type metrics perform well for hyperbolic PDEs, consistent with previous tests,
showing strong correlation with simulation error in Table~\ref{spearmanHyper}.
In particular, IMPF and RKPF, which already ranked among the best performing quality metrics in earlier studies,
again emerge as the leading quality metrics within this type, emphasizing  their robustness and reliability across different PDE types.

%
%
%
%

Within their respective categories, DJM (\textit{mesh} type), FD  (\textit{geometry} type),
and IME (\textit{energy} type) exhibit the highest monotonic dependency on simulation error. These results
also align with the findings from the elliptic PDE study above.

Empirically, we observed that the metric correlations were stronger with the root mean square error in this case.
However, we only report the infinity norm error for consistency with the elliptic problem tests.


\subsection{Discussion}
\label{sec:Discussion}

In Sections~\ref{parameterStudy} to~\ref{sec:HyperbolicResults}, we conducted extensive numerical testing
to evaluate the correlation between various quality metrics and simulation error across multiple scenarios.
We now want to identify the ``best" quality metrics as those that are the most robust and reliable, exhibiting a
strong correlation with numerical error across different test cases and scenarios.
Fortunately, this task is simplified by the fact that the performance score $\sigma$ enables a
quantitative comparison across different settings. In particular, due to its remarkable stability  in
all experiments for the Poisson equation in 2D and 3D, we present in Table~\ref{scores} a single ``elliptic"
score for each metric.  
Due to differences in the performance of a few metrics, such as  ST and DW,  between elliptic and
hyperbolic cases, ``hyperbolic"  scores, copied from Table~\ref{spearmanHyper}, are provided in a separate column
of  Table~\ref{scores} for the ease of comparison.

\begin{table}[h!]
	\centering
		\caption{Performance score of the quality metrics.
		The column \emph{Elliptic} shows the average of $\bar\sigma$
		over all Poisson experiments
		(Tables~\ref{tab:spearma-setsize-u1}--\ref{tab:3d-u1} and
		\ref{tab:spearma-subsetsize-u2}--\ref{tab:spearma-3d-u2}),
		whereas
		\emph{Hyperbolic} gives $\sigma$ of Table~\ref{spearmanHyper}.
		}
	\label{scores}
	\begin{tabular}{lcccc}
	\toprule
		\textbf{Metric}  & \textbf{Type} & \textbf{Elliptic} & \textbf{Hyperbolic} \\
	\midrule
 \multicolumn{4}{l}{\textbf{Energy}} \\
RKE & energy   &  5.88 &  4 \\
LKE  & energy   & 29.00 & 40 \\
IME  & energy   &  6.00 &  4 \\
WKE  & energy   & 21.50 & 26 \\
CE   & energy   & 31.25 & 17 \\
\midrule
\multicolumn{4}{l}{\textbf{Equilibrium}} \\
RKPF & equilibrium &  5.50 &  5 \\
LKPF & equilibrium & 12.38 &  9 \\
IMPF & equilibrium &  3.13 &  4 \\
WKPF & equilibrium & 12.13 &  10 \\
CPF  & equilibrium & 13.38 & 11 \\
\midrule
\multicolumn{4}{l}{\textbf{Geometry}} \\
SD   & geometry & 27.75 & 46 \\
FD   & geometry &  2.25 &  3 \\
UR   & geometry & 26.75 & 52 \\
PD   & geometry & 29.25 & 42 \\
ALR  & geometry & 25.75 & 16 \\
RLR  & geometry &  9.13 &  8 \\
\midrule
\multicolumn{4}{l}{\textbf{Mesh}} \\
IAR   & mesh     & 16.00 & 23 \\
DJM  & mesh     &  2.63 &  1 \\
MCD  & mesh     & 10.50 &  8 \\
\midrule
\multicolumn{4}{l}{\textbf{System}} \\
GFL1 & system   &  2.00 &  2 \\
GFL2 & system   &  1.88 &  2 \\
CN   & system   & 18.00 & 37 \\
RE   & system   &  1.38 &  2 \\
ST   & system   & 66.88 & 128 \\
DW   & system   & 11.38 & 159 \\
\bottomrule
		\end{tabular}
	\end{table}

Based on Table~\ref{scores} the metrics with the \emph{best overall performance} (both scores not exceeding 4) 
are RE, GFL2, GFL1, DJM, FD, and IMPF
 summarized in Table~\ref{Hperform}.
Indeed, they consistently prove to be the best error indicators, and thus, good indicators for point cloud quality in all
experiments, for both elliptic and hyperbolic problems.

\begin{table}[!htp]
\caption{Summary of the best performing point cloud quality metrics.
	\emph{Description} refers to the section in Appendix~\ref{qualityMetricDetails} where each metric is explained in detail.
	\emph{Maps} indicate the best  transformation from local to global for each quality metric.
	The last column indicates the computational cost of each metric.}
	\centering
		\begin{tabular}{lcccc}
			\toprule
			Metric & Maps &Type & Description &	Cost \\
			\midrule
		Recovery Error (RE) & std, rms & System Matrix & \ref{maxAppErro}        & Moderate \\
		Growth Function via $\ell_2$-norm (GFL2) & std, rms, mean & System Matrix & \ref{growthF} & Moderate \\
		Growth Function via $\ell_1$-norm (GFL1) & std, rms, mean & System Matrix & \ref{growthF} & Moderate \\
		Determinant of the Jacobian Matrix (DJM) & std            & Mesh          & \ref{jacobM} & High \\
		Fill Distance (FD)                       & std, rms       & Geometry      & \ref{fillDistance} & High \\
		Inter-Molecular Potential Force (IMPF)   & mean, rms      & Equilibrium   & \ref{potentialForceIMP} & Low \\
		

			\bottomrule
		\end{tabular}%
	\label{Hperform}
\end{table}

For these metrics to be used in practice, to improve point cloud quality, for point cloud generation,
refinement, or adaptation after a Lagrangian movement process, their computational efficiency must also be considered.
In this sense, IMPF is an unbeatable leader because its
computational cost is minimal as it is computed directly from
the locations of nearby points without solving any linear systems.
On the other hand, the metrics RE, GFL1, and GFL2 require higher computational effort because they 
require computing numerical differentiation weights of Equation~\eqref{diffFormula} by solving small linear systems.
Note that the differentiation weights form the system matrix of \eqref{discritPoisson}, and therefore they may also 
be needed for other tasks of the meshfree workflow.
The cost of both DJM and FD is even significantly higher because they
require a mesh to be computed over the point cloud. (See Appendix~\ref{qualityMetricDetails}, in particular 
Section~\ref{fillDistance} for the geometry-based metric FD, whose most efficient implementation 
uses a Delaunay triangulation, making it a mesh-based metric in this sense.)
Even though computing local meshes can speed up this process, particularly when executed in parallel,
this introduces undesired computational overhead, and might require expensive mesh improvement in some situations.
This difference in the computational cost  of the metrics is indicated in the last column of Table~\ref{Hperform}.

We also identify several metrics with a \emph{good overall performance}, namely  RKE, IME and RKPF, whose
performance scores in Table~\ref{scores} are between 5 and 6 in the elliptic case and between 4 and 5 in the 
hyperbolic case. They all are of a low computational cost.

The next group includes metrics RLR, MCD, WKPF, LKPF and CPF, with the performance score between 9 and 13 in the elliptic case
and  between 8 and 11 in the hyperbolic case, which we evaluate as an \emph{acceptable overall performance}.
Among them RLR, WKPF, LKPF and CPF are of low computational cost, but MCD requires a mesh.

We now highlight differences in performance within the
\emph{types} of metrics introduced in Section~\ref{typesAndMetrics}, with some unexpected results.
Note that every type includes both good and poor performers.

\begin{itemize}
\setlength{\itemsep}{1pt}
\setlength{\parskip}{0pt}
\setlength{\parsep}{0pt}


\item  Of the \textit{energy} type, RKE and IME perform well, although they do not belong to the very best metrics.
We had expected a similar performance from RKE and LKE, given the well known and widely used 
Riesz and Log kernel energies. However, the results for LKE are disappointing, which makes it the 
worst performing metric within this type.
%
%


	\item The results suggest that the \textit{equilibrium} type quality metrics generally exhibit a stronger monotonic 
	correlation with the error than their corresponding \textit{energy} type counterparts. 
	Among them, IMPF ranks among the best-performing metrics overall, while RKPF belongs to the good  performers, 
	and the remaining	three metrics all show an acceptable result. See also a comment below on WKPF that significantly improves
	its performance when computed on a reduced support.



\item Among the \textit{geometry}-based quality metrics, FD is clearly the best performer and ranks among the very best
metrics overall, while the results for RLR are acceptable. 
In contrast, ALR, SD and PD, although well established in point cloud quality assessment and point cloud generation, 
show a rather weak performance. 
%

    \item Within the \textit{mesh} based quality metrics, DJM emerges as a reliable metric for both elliptic and hyperbolic PDEs, and MCD shows acceptable
		performance. 
		While aspect ratio IAR is a well-established and very reliable metric to quantify mesh quality, we observe that
	it does not serve as a good predictor of point cloud quality, in particular in all 2D tests. 
	This highlights the fundamental difference between quality assessment of a point cloud and a mesh,  
	and shows that point cloud quality assessment would not be 
	well served by meshing the point cloud and then assessing the quality of the resulting mesh. 
	Furthermore, as stated earlier, \textit{mesh}-type metrics are computationally expensive because they require 
	the construction of a mesh, which is not otherwise needed in the mesh-free methods.

%

\item The \textit{system-matrix} metrics RE, GFL1, and GFL2 rank among the very best performers in both the elliptic and
hyperbolic tests, making them the most reliable metrics within this category. We expected that they would significantly
outperform other metrics for higher order
methods because they estimate the error of the  numerical differentiation weights that change significantly with the
polynomial order and are  not available to other metrics. Interestingly, the very simple system-matrix metric DW that only
measures the diagonal weights performs very well in all low order tests for the Poisson equation, with average
$\bar\sigma=4.17$ over these experiments. However, it performs weakly for higher order methods
(Tables~\ref{tab:spearma-porder-u1} and 	\ref{tab:spearma-porder-u2}) and very poorly for the hyperbolic problem.
	In our preliminary experiments the stability metric ST correlated well  with the norm of the inverse of  the system matrix
(the classical stability constant for the error analysis of the Poisson problem), but its correlation with
the numerical error is quite poor.
\end{itemize}
%

The best aggregation maps to obtain global metrics from the local metrics are \texttt{rms},
\texttt{mean}, and \texttt{std}. This is in contrast with mesh quality assessment, where extreme
values of local metrics (\texttt{max} or \texttt{min}) are the most widely used.
In the \emph{Maps} column of Table~\ref{Hperform}  we included recommended aggregation maps for the best performing metrics.
Notably, for RE, GFL2 and some other point quality metrics more than one aggregation map can be used with almost
the same results. For example, the correlation coefficients $\rho$ of the global metrics obtained with \texttt{rms} and
\texttt{mean} with the error remains nearly unchanged for IMPF, with variations of $\sigma$ between the maps not exceeding one.
%

Several variations of the metrics considered here are possible. Different aggregation maps can be used to obtain the global metric from local values. 
For some metrics, such as DJM and IAR, there is an additional choice in how the local value itself is computed from the simplices considered. 
Another possible variation is to evaluate a metric only on a subset of the support, rather than on the full support.
In this work, we compute all local quality metrics using the full support \(\mathbf{S}_i\), in order to keep the comparison uniform. 
However, we observed that the metrics WKE and WKPF based on a smoothing kernel, can show stronger correlation with numerical error when 
computed on a reduced support. This indicates that these metrics may incorrectly account for points   near boundaries of
supports. 
Restricting to the interior part of the support significantly improves their performance, 
apparently by reducing these boundary effects. In particular, the values of $\bar{\rho}$ for WKPF computed with half of the support size
increase to 0.960 and 0.859 for Tables~ \ref{tab:3d-u1} and \ref{spearmanHyper}, respectively. This results in a significant
improvement in the performance score $\bar{\sigma}$, which decreases from 16 to 3 for Table~\ref{tab:3d-u1} and from 10 to 5
for Table~\ref{spearmanHyper}.
It may also be worthwhile to investigate whether kernel correction techniques \cite{bonet1999variational} 
can improve the performance of 
these metrics when they are evaluated on the full support. 
A detailed study of these variants would go beyond the scope of the present manuscript. 
Here, our aim was to provide a consistent baseline for comparing a broad range of point cloud quality metrics.
To this end, we already consider several variations in the numerical setting, 
including changes in resolution, support size, 
polynomial order, spatial dimension, and PDE type. 


We emphasize that all simulations conducted in this work, and thus, all conclusions reached, are based
on a specific class of meshfree collocation methods, namely, meshfree Generalized Finite Difference Methods (GFDMs).
It is possible that for other meshfree methods the best-performing metrics may be different from those
identified in this study, or  some of our conclusions may not hold in different meshfree frameworks.

 	\section{Conclusion}
\label{sec:Conclusion}

Ensuring good quality point cloud discretizations is essential for achieving accurate
and stable numerical solutions in meshfree methods. However, unlike mesh-based methods,
where well-established metrics exist to quantify mesh quality, a systematic study of point
cloud quality has been lacking. In this work, we addressed this gap by thoroughly
investigating a wide range of quality metrics, covering existing ones from the literature
as well as several new metrics introduced in this study. We considered 25 local quality
metrics and applied 7 different aggregation maps to each of them to obtain global metrics. Given
the large number of metrics, we introduced a classification system to group similar metrics
into five categories: geometry-based, energy-based, equilibrium-based, mesh-based, and
system matrix-based metrics. This classification provides a structured way to investigate,
compare and analyze different quality indicators.

To determine which of these metrics best quantify point cloud quality, we proposed that a
reliable metric should consistently predict numerical error across diverse scenarios. Focusing
on meshfree collocation methods, we evaluated the performance of these metrics in a variety of
settings: for both elliptic and hyperbolic PDEs, in 2D and 3D domains, across different
resolutions, support sizes, and polynomial orders.

Our approach to the assessment of the quality metrics included testing many randomly generated point clouds for
each scenario, measuring monotonic correlation between a metric and the ground truth error of a known analytical solution, 
and comparing the ranks of the metrics according to the values of the correlation coefficients for each scenario.

Our numerical study revealed that the ranks (performance scores) of the well-behaving metrics are remarkably stable between
various tests, which allowed us to identify six metrics that consistently serve as robust indicators of
numerical accuracy, and thus point cloud quality, across all scenarios considered:
\begin{itemize}
\setlength{\itemsep}{1pt}
\setlength{\parskip}{0pt}
\setlength{\parsep}{0pt}
	\item Recovery Error (RE)
\item Growth Function via $\ell_2$-norm (GFL2)
\item Growth Function via $\ell_1$-norm (GFL1)
\item Determinant of the Jacobian Matrix (DJM)
\item Fill Distance (FD)   
\item 	Inter-Molecular Potential Force (IMPF)
\end{itemize}


An important observation was that several other metrics
commonly used in the literature to assess point cloud quality show a weak correlation with the numerical error.

Through this systematic analysis, we provided a comprehensive foundation for answering the
question: \textit{what constitutes a good point cloud?} Our findings can guide the development
of adaptive algorithms for point cloud generation and refinement, thus improving the accuracy
and efficiency of meshfree simulations. 

 	\appendix
 	\begin{appendices}

\section{Point cloud quality metrics}\label{qualityMetricDetails}

In this appendix, we give detailed descriptions of all the \emph{local} metrics
considered for evaluating point cloud quality. See Section~\ref{metrics} and
Table~\ref{tab:qualityMetrics} for a list of all the metrics, and for the aggregation maps to
define the global quality metrics from local ones. In Section~\ref{typesAndMetrics},
we introduced five types of quality metrics to group together similar definitions and
allow a clearer understanding of their characteristics.
We emphasize that this classification is based on the metric definitions, rather than on  the
procedure to compute the metric.
Here, we also outline the motivation for each type of metric, and then present the definitions
and mathematical formulations of the metrics. Implementation details are included only when
the computations are not straightforward.
The metrics are presented here in the same order as that in Table~\ref{tab:qualityMetrics}.

We emphasize that all metrics considered here are locally defined. To recall the notation used:
For a point $\mathbf{x}_i$ in a point cloud $\mathbf{X} \subset \mathbb{R}^\nu$, all metrics are computed over the
support or neighborhood $\mathbf{S}_i$ consisting of $m=m(i)$ nearby points, where $m$ can vary
from point to point. $\mathbf{S}_i$ is determined by proximity,
see Equation~\ref{range}. Note that some metrics use all the points in the
neighborhood, while others are dependent only on a subset of $\mathbf{S}_i$.

\subsection{Energy-based quality metrics}\label{energyType}

This type of metrics treat the point cloud as a system of interacting particles. Point cloud quality
is quantified as the total system energy, based on pair-wise particle-particle interactions.
The assumption is that lower energy typically corresponds to more uniform and well-spaced point distributions.

Given an interaction potential $\Psi = \Psi(\mathbf{x}_i,\mathbf{x}_j)$, defined between a pair of points $i$ and $j$,
an energy based quality metric $E$ at point $\mathbf{x}_i$ is defined as the total interaction energy
\begin{equation}
\label{Eq:EnergyTypeMetric}
	E(\mathbf{x}_i) =
	\sum_{\substack{\mathbf{x}_j \in \mathbf{S}_i, \\ \mathbf{x}_j \neq \mathbf{x}_i}}
	\Psi(\mathbf{x}_i, \mathbf{x}_j) .
\end{equation}

Different metrics of this type consider different functions for the interaction potential $\Psi$.
The first three metrics of this type, namely the Riesz kernel energy, logarithmic kernel energy,
and inter-molecular energy, have been used as a part of point cloud generation algorithms in meshfree literature.
These methods begin with a point cloud generated using an unrelated approach of the user's choice (see~\cite{suchde2023point}),
and then in a second step, move the points to minimize the chosen energy function. This is based on the assumption that the
lower energy corresponds to a ``better'' point cloud.
In this work, we repurpose the energy functions used in these point cloud generation methods as
stand-alone metrics to directly quantify point cloud quality, independent of any generation process.
Based on similar ideas, we also introduce two additional energy functions to serve as point cloud quality
indicators: the Wendland energy and the centroid energy.

\subsubsection{Riesz kernel energy (RKE)}\label{rieszEnergy}

The Riesz kernel energy (RKE) is a generalization of the Coulomb electrostatic potential. For the RKE, the
potential in Equation~\eqref{Eq:EnergyTypeMetric} is taken as $\Psi = \|\mathbf{x}_i - \mathbf{x}_j\|^{-s} $,
for an integer parameter $s$. Thus, we have the local metric
\begin{equation}
	\text{RKE}(\mathbf{x}_i) =
	\sum_{\substack{\mathbf{x}_j \in \mathbf{S}_i, \\ \mathbf{x}_j \neq \mathbf{x}_i}}
	\|\mathbf{x}_i - \mathbf{x}_j\|^{-s}.
	\label{reiszEnergy}
\end{equation}

This energy function has been used, for example, in \cite{vlasiuk2018fast} to generate quasi-uniform point
distributions. Various values of the parameter $s$ have been considered in literature.
It has been shown \cite{landkof1972foundations} 
that point clouds
that minimize the RKE are asymptotically uniform, in the sense that as the number of points increases,
the points tend towards being uniformly distributed. This property is helpful in the theoretical
convergence analysis of meshfree methods.


We tested several values of $s$ in our numerical experiments in the range $0 < s \leq \nu +1$.
We observed that $s$ has a significant effect on the correlation between the RKE metric and the numerical error. 
Among all cases considered, $s=3$ in 2D and $s=4$ in 3D gave the strongest correlations. 
Thus, all RKE results reported in this work use $s=3$ for 2D domains and $s=4$ for 3D domains.


\subsubsection{Log-kernel energy (LKE)}\label{logaEnergy}
Based on motivations in potential theory and the study of harmonic functions \cite{landkof1972foundations},
the Log-kernel energy (LKE) is defined as
\begin{equation}
	\text{LKE}(\mathbf{x}_i) =
	\sum_{\substack{\mathbf{x}_j \in \mathbf{S}_i, \\ \mathbf{x}_j \neq \mathbf{x}_i}}
	\|\mathbf{x}_i -\mathbf{x}_j\|
	\left( \log \frac{ \|\mathbf{x}_i - \mathbf{x}_j\|}{2} - 1 \right) + 2.
	\label{logEnergy}
\end{equation}
%

This energy function has been used in \cite{kunc2019generation} to generate quasi-uniform point
distributions. Similar to the case of RKE, it has also been shown that sequences of point clouds that
minimize the LKE are asymptotically uniform \cite{kunc2019generation}.

\subsubsection{Inter-molecular energy (IME)}\label{IMenergy}

The Lennard-Jones potential (LJP) is one of the most widely studied models for inter-molecular (IM) interactions.
The LJP has been used as an energy function in node placement for triangular mesh generation, for example,
by \cite{zhang2005node}. Inspired by the LJP, similar IM interactions have also been used for
meshfree point cloud generation~\cite{negi2021algorithms}.
Based on \cite{negi2021algorithms}, we define the local metric inter-molecular energy
for a point \(\mathbf{x}_i\) as
\begin{equation}
	\text{IME}(\mathbf{x}_i) =
	12 K_r\sum\limits_{\substack{\mathbf{x}_j \in \mathbf{S}_i, \\ \mathbf{x}_j \neq \mathbf{x}_i}}
	\left( \frac{c^2}{r_{ij}^3} - \frac{c}{r_{ij}^2} \right),
	\qquad \text{ with } r_{ij} := \|\mathbf{x}_i - \mathbf{x}_j\|,\quad c = \tfrac23\alpha \Delta s,
	\label{Eq:IMenergy}
\end{equation}
%
%
%
where \(\Delta s\) depends on the particle spacing,
 \(\alpha = 0.95\) is a scaling factor, and $K_r = 0.004 \Delta s$  following \cite{negi2021algorithms}.  
Furthermore, \(r_{ij}\) denotes the Euclidean distance between the points \(\mathbf{x}_i\) and
\(\mathbf{x}_j\).
We tested different values of \(\Delta s\), before setting $\Delta s = \lambda$, which exhibited the highest correlation of metric with error.



\subsubsection{Wendland kernel energy (WKE)}\label{WandlandEnergy}

In many meshfree and particle-based methods, function values are approximated using weighted
sums over neighboring points. A minimum requirement for these approximations to be accurate is that the weights
assigned to neighboring points should sum to unity. Uniformly distributed point clouds typically
satisfy this requirement, with deviations increasing for irregular point distributions, see, for example,
\cite{liu2003smoothed}.
Motivated by this idea, we propose an energy-based quality metric based on the quintic Wendland kernel
\cite{wendland2004scattered}, 
which is a common weighting function in particle methods (for example, \cite{abdolahzadeh2019mixing}).
\begin{equation}
	\text{WKE}(\mathbf{x}_i) =
	\sum_{\substack{\mathbf{x}_j \in \mathbf{S}_i \\  \mathbf{x}_j \ne \mathbf{x}_i}}
	\alpha_\lambda\left( 1-\frac{q}{2} \right)^{4}\left(2q+1\right),\qquad q = \frac{r_{ij}}{\lambda},
	\label{WandlandPotentialEnergy}
 \end{equation}
where $\lambda$ is  the support size/interaction radius \eqref{range}, and \(r_{ij}\) is the Euclidean distance
between the particles. Parameter $\alpha_\lambda$ depends on $\lambda$ and the spatial dimension of the domain,
\(\alpha_\lambda\) is \(\frac{7}{4\pi \lambda^2}\) in 2D and \(\frac{21}{16\pi \lambda^3}\) in 3D.

Based on the motivation in function approximation in particle methods, we can say that
WKE serves as a quality metric to evaluate the regularity of point distributions.

\subsubsection{Centroid energy (CE)}\label{balaceKernelEnergy}

We also propose another energy-based metric, the centroid energy (CE), defined as
\begin{equation}
	\text{CE}(\mathbf{x}_i) = \frac{1}{m} \sum_{\substack{\mathbf{x}_j \in \mathbf{S}_i \\ \mathbf{x}_j \ne \mathbf{x}_i}}
\frac{r_{ij}^2}{2},
	\label{balanceEnergy}
\end{equation}
where $ m=|\mathbf{S}_i|$ is the number of points in the support and $r_{ij}$ is the Euclidean distance between points.

\subsection{Equilibrium-based quality metrics}

Metrics of this type are derived from the gradient of potential functions used in the energy-based
metrics. They measure how close a point cloud is to a force-balanced or equilibrium state.
For an interaction potential $\Psi$ used to define an energy-based quality metric, the corresponding potential force
$\mathbf{F}_{\Psi}$ at a point \(\mathbf{x}_i\) is given by
\begin{equation}
	\mathbf{F}_{\Psi}(\mathbf{x}_i) =
	\sum_{\substack{\mathbf{x}_j \in \mathbf{S}_i \\ \mathbf{x}_j \ne \mathbf{x}_i}}
	\nabla \Psi(\mathbf{x}_i,\mathbf{x}_j)
\label{equilibriumType}
\end{equation}
where \( \nabla \Psi \) denotes the gradient of \(\Psi\).
The equilibrium-type quality metric is then defined as \( \|\mathbf{F}_{\Psi}(\mathbf{x}_i) \| \).
This can be interpreted as the norm of the net potential force acting on point \(\mathbf{x}_i\). As a point moves
closer to equilibrium, the net force on it will decrease. Low values of any of these metrics can be considered as
a point cloud close to equilibrium, thus these metrics serve as a potential indicator for point cloud quality.

This type of metric was motivated by the use of the inter-molecular energy for point cloud generation in energy
form (Appendix~\ref{IMenergy}) in \cite{zhang2005node}, and in the equilibrium form in
\cite{smirnov2008physically, negi2021algorithms}.
We extend this to obtain an equilibrium-based metric for each of the energy-based metrics considered above
 (see Appendix~\ref{energyType}).
It is important to note that while the corresponding energy-type and equilibrium-type metrics are related
and can be derived from each other, they are not the same.

\subsubsection{Riesz kernel potential force (RKPF)}\label{rieszPF}

Taking the gradient of Equation~\eqref{reiszEnergy}, we get the force
\begin{equation}  \label{RKPF}
	\mathbf{F}_{\text{RKE}}(\mathbf{x}_i) =
	\sum_{\substack{\mathbf{x}_j \in \mathbf{S}_i \\ \mathbf{x}_j \ne \mathbf{x}_i}}
	 -s  \|\mathbf{x}_i - \mathbf{x}_j\|^{-(s+1)} \mathbf{e}_{r_{ij}}, \qquad \text{with}\quad
	\mathbf{e}_{r_{ij}} := \frac{\mathbf{x}_i - \mathbf{x}_j}{\|\mathbf{x}_i - \mathbf{x}_j\|},
\end{equation}
where \(\mathbf{e}_{r_{ij}}\) is the unit vector in the direction from \(\mathbf{x}_i\) to \(\mathbf{x}_j\).
The corresponding equilibrium-type quality metric, the Riesz kernel potential force (RKPF), is then defined as
\begin{equation} \label{MRKPF}
	\text{RKPF}(\mathbf{x}_i) = \|\mathbf{F}_{\text{RKE}}(\mathbf{x}_i)\|.
\end{equation}
Similar to the RKE metric, we tested different values of $s$ and found that $s=1$ in 2D and $s=4$ in 3D give the strongest 
correlation between the RKPF metric and the numerical error. Thus, all RKPF results reported in 
this work use $s=1$ for 2D domains and $s=4$ for the 3D domains.

\subsubsection{Logarithmic kernel potential force (LKPF)}\label{logPF}

Taking the gradient of Equation~\eqref{logEnergy}, we get the force
\begin{equation}\label{LKPF}
	\mathbf{F}_{\text{LKE}}(\mathbf{x}_i) = - \sum_{\substack{\mathbf{x}_j \in \mathbf{S}_i \\ \mathbf{x}_j \ne
	\mathbf{x}_i}}  \log \frac{\|\mathbf{x}_i - \mathbf{x}_j\|}{2}   \mathbf{e}_{r_{ij}}.
\end{equation}
The corresponding equilibrium-type quality metric, the Logarithmic kernel potential force (LKPF), is then defined as
\begin{equation} \label{MLKPF}
	\text{LKPF}(\mathbf{x}_i) = \|\mathbf{F}_{\text{LKE}}(\mathbf{x}_i)\|.
\end{equation}

\subsubsection{Inter-molecular potential force (IMPF)}\label{potentialForceIMP}

The inter-molecular energy IME (Appendix~\ref{IMenergy}) has also been used for point cloud generation
in a force form in \cite{smirnov2008physically, negi2021algorithms}. Specifically, \cite{negi2021algorithms}
considers the force
%
\begin{equation}\label{IMPF}
	\mathbf{F}_{\mathrm{IME}}(\mathbf{x}_i)=
	\begin{cases}
		\sum\limits_{\substack{\mathbf{x}_j \in \mathbf{S}_i \\ \mathbf{x}_j \ne \mathbf{x}_i}}
		192\,K_r\left(\frac{3c^2}{\Delta s^4}-\frac{c}{\Delta s^3}\right)\mathbf{e}_{r_{ij}},
		& r_{ij}\le \dfrac{\Delta s}{2}, \\[1.2ex]
		\sum\limits_{\substack{\mathbf{x}_j \in \mathbf{S}_i \\ \mathbf{x}_j \ne \mathbf{x}_i}}
		12\,K_r\left(\frac{3c^2}{r_{ij}^4}-\frac{2c}{r_{ij}^3}\right)\mathbf{e}_{r_{ij}},
		& \dfrac{\Delta s}{2}<r_{ij}<\alpha\Delta s, \\[1.2ex]
		\mathbf{0},
		& r_{ij}\ge \alpha\Delta s.
	\end{cases}
\end{equation}
where $\Delta s$ is set equal to $\lambda$, consistent with the definition used for the IME metric (Appendix \ref{IMenergy}). Furthermore,  \(c = \frac{2\alpha \Delta s}{3}\),
with \(\alpha = 0.95\) as a scaling factor. Additionally,
\(r_{ij}\)
denotes the distance between particles \(\mathbf{x}_i\) and \(\mathbf{x}_j\), and \(\mathbf{e}_{r_{ij}}\) is the unit
vector in the direction from \(\mathbf{x}_i\) to \(\mathbf{x}_j\). All these constants and parameters are taken as in
Appendix \ref{IMenergy}. 

We note that unlike the case of other equilibrium-based quality metrics defined in this section,
$\mathbf{F}_{\text{IME}} \neq \nabla \text{IME}$. The formulation of $\mathbf{F}_{\text{IME}}$ considered here is based on that considered
in~\cite{negi2021algorithms}, which uses the gradient derived in an SPH discretization.

Based on this force, we define the inter-molecular potential force (IMPF) quality metric as
\begin{equation} \label{MIMPF}
	\text{IMPF}(\mathbf{x}_i) = \|\mathbf{F}_{\text{IME}}(\mathbf{x}_i)\|.
\end{equation}

\subsubsection{Wendland kernel potential force (WKPF)}\label{wandPF}

Taking the gradient of Equation~\eqref{WandlandPotentialEnergy}, we get the force
\begin{equation}
	\mathbf{F}_{\text{WKE}}(\mathbf{x}_i) =
	\frac{- 5 \alpha_\lambda }{ \lambda^2}
	\sum_{\substack{\mathbf{x}_j \in \mathbf{S}_i, \\ \mathbf{x}_j \neq \mathbf{x}_i}}
	 \left( 1 - 0.5 \frac{r_{ij}}{\lambda} \right)^3 (\mathbf{x}_i - \mathbf{x}_j).
	\label{potentialForceWandland}
\end{equation}
We introduce the Wendland Kernel Potential Force (WKPF) local quality metric, defined as:
\begin{equation}
	\text{WKPF}(\mathbf{x}_i) = \|\mathbf{F}_{\text{WKE}}(\mathbf{x}_i) \|.
\end{equation}
%

\subsubsection{Centroid potential force (CPF)}\label{baPF}

All the forces considered above in the equilibrium-based metrics can be interpreted to move the
point \( \mathbf{x}_i \) toward a weighted centroid of its neighbors.
Based on this, we introduce another force that simply pulls \( \mathbf{x}_i \) toward the
centroid, without any additional weights. For this, we introduce the centroid force
by taking the gradient of Equation~\ref{balanceEnergy}
\begin{equation}
	\mathbf{F}_{\text{CE}}(\mathbf{x}_i) =
	\frac{1}{m}
	\sum_{\substack{\mathbf{x}_j \in \mathbf{S}_i \\ \mathbf{x}_j \ne \mathbf{x}_i}}
	(\mathbf{x}_i - \mathbf{x}_j)  ,
	\label{CPFForce}
\end{equation}
where $m=|\mathbf{S}_i|$ is the number of points in the support, which serves as a
normalization factor.
The corresponding metric, the centroid potential force (CPF), is defined as
\begin{equation}
	\text{CPF}(\mathbf{x}_i) =  \|\mathbf{F}_{\text{CE}}(\mathbf{x}_i) \|.
\end{equation}
In this case the  equilibrium condition is that the geometric centroid of all neighboring points
coincides with the center point $\mathbf{x}_i$.

\subsection{Geometry-based quality metrics}

These metrics are computed directly from the relative locations of points in the support domain,
without relying on any intermediate constructs such as energy functions, force balances, or meshes.
They are based purely on spatial properties like pairwise distances, nearest neighbors, or local
density. The separation distance and fill distance metrics in this type are among the most widely
used point cloud quality metrics, and are also commonly used in theoretical studies for proving
convergence and error bounds.

\subsubsection{Separation distance (SD)}\label{Separation}

The separation distance (SD) for a point $\mathbf{x}_i$ with support $\mathbf{S}_i$ is defined as
\begin{equation}
	\text{SD}(\mathbf{x}_i) = \frac{1}{2}
	\min_{ \substack{ \mathbf{x}_j,\mathbf{x}_k\in \mathbf{S}_i \\ \mathbf{x}_j \ne \mathbf{x}_k } }
 	\| {\mathbf{x}_j - \mathbf{x}_{k}} \|.
	\label{kerApprox}
\end{equation}

The separation distance is half the minimum distance between any two distinct points in $\mathbf{S}_i$.
Imposing a lower bound on the separation distance is a widely used strategy for ensuring good-quality point clouds~\cite{drumm2008finite, suchde2019fully},
as it prevents points from clustering too closely. The separation distance also plays a crucial role in the
theoretical analysis of many meshfree methods \cite{fasshauer2007meshfree}.
In particular, it has been demonstrated that a small separation
distance adversely affects the condition number of the radial basis function interpolation matrix,
which is a standard criterion for assessing numerical stability \cite{fasshauer2007meshfree}.
Therefore, SD is a natural choice for a metric to assess point cloud quality.
It has also been used across meshfree literature to analyze point cloud regularity
\cite{duh2021fast, duh2024discretization, fornberg2015fast, shankar2018robust}.

\subsubsection{Fill distance (FD)}\label{fillDistance}

The fill distance (FD) for a point $\mathbf{x}_i$ with support $\mathbf{S}_i$ is defined as
\begin{equation}
	\text{FD}(\mathbf{x}_i) =
	\adjustlimits
	\sup_{\mathbf{x} \in B_i} \min_{\mathbf{x}_j \in \mathbf{S}_i}
	\|\mathbf{x} - \mathbf{x}_j\|,
\end{equation}
where $B_i$ is the smallest ball centered at $\mathbf{x}_i$ that contains all points of $\mathbf{S}_i$.
Thus, the fill distance is the radius of the largest ball  centered in the neighborhood of $\mathbf{x}_i$
that does not contain any point from $\mathbf{S}_i$ \cite{wendland2004scattered}.

Similar to the separation distance, the fill distance is also widely used across meshfree literature.
It has been used as a quality metric \cite{duh2024discretization}, and
enforcing an upper bound on the fill distance has been used to ensure good quality point clouds in practice
\cite{drumm2008finite, suchde2019meshfree}. It has also been used extensively in theoretical analysis, for example
in error estimation \cite{fasshauer2007meshfree, schaback2015computational,wendland2004scattered}.

Unlike the metrics considered so far, computation of FD is not trivial. Various approaches to the
computation or estimation of FD have been proposed, including both algebraic approaches based on a Euclidean
distance matrix \cite{fasshauer2007meshfree}, and geometric approaches, that have also been referred to as
hole detection \cite{suchde2019fully}. Among the geometric approaches, some use Delaunay
triangulations \cite{du2002meshfree, suchde2019fully}, others their dual Voronoi diagrams \cite{drumm2008finite},
while some other approaches use discrete hole searching methods \cite{suchde2023point} that do not require any mesh.


Throughout this work, the fill distance is computed geometrically using a Delaunay triangulation of the point cloud.
For each point $\mathbf{x}_i$, we collect all Delaunay simplices within its support $\mathbf{S}_i$ and compute their  circumradii. 
The local fill distance is then estimated as the largest of these circumradii among these simplices.
We also tested other approaches to compute the FD metric, which produced very similar results.

We note that despite using a mesh to compute FD, we classify it under a geometry-based quality metric, and
not a mesh-based one, since the definition of FD does not require a mesh, and hence its performance is not supposed to depend
on the mesh quality.
On the other hand, we are not aware of a cheap meshfree approach to compute
or estimate FD that outperforms the mesh-based ones. Note that the Delaunay
triangulations or Voronoi tessellations can be computed either locally on $\mathbf{S}_i$ or globally on
the entire point cloud.

\subsubsection{Uniformity ratio (UR)}\label{meshRatio}


The ratio of fill distance to separation distance appears in many theoretical convergence and stability
proofs and error analysis in meshless literature \cite{fasshauer2007meshfree,wendland2004scattered}.
Here, we also use it as a local metric to assess point cloud quality,
\begin{equation}
	\text{UR}(\mathbf{x}_i) = \frac{\text{FD}(\mathbf{x}_i)}{\text{SD}(\mathbf{x}_i)}.
\end{equation}
This ratio has also been referred to as the mesh ratio and mesh separation ratio in meshfree literature.
Here, we adopt the term \emph{uniformity ratio} to avoid any confusion suggesting that a mesh is required.
This ratio has also been used to judge point cloud quality \cite{van2021fast}.

A smaller value of UR indicates that the fill distance and separation distance are comparable,
meaning points are evenly spaced without excessive clustering or large gaps.
In the asymptotic theoretical estimates often the existence of an upper bound on UR is assumed, a property of point clouds
referred to as quasi-uniformity.


\subsubsection{Packing density (PD)}\label{SecPackDensity}

Another way to quantify point cloud quality is to consider how densely the points are packed in
the domain. This idea, rooted in classical sphere or circle packing problems \cite{conway2013sphere}, evaluates
how efficiently the space is filled while maintaining a minimum distance between points. Based on the use of
packing density as a quality indicator in literature \cite{van2021fast}, we define the local packing density
$\eta_i$ at point $\mathbf{x}_i$ as
\begin{align}
	\eta_i &= \frac{m}{\left|B_i\right|}
	\frac{\pi^{\nu/2} }{\Gamma\left(\frac{\nu}{2} + 1 \right)}
	\text{SD}(\mathbf{x}_i)^\nu,  \\
	&= 	m \left( \frac{\text{SD}(\mathbf{x}_i)}{R_{B_i}} \right)^{\nu},
	\label{pointSetPaking}
\end{align}
where $|B_i|$ is the volume of the smallest ball $B_i$ centered at $\mathbf{x}_i$ that contains all points of $\mathbf{S}_i$,
SD is the separation distance (see Appendix~\ref{Separation}),
$m=\left|\mathbf{S}_i\right|$ is the number of points in the support, and $\Gamma$ is the Gamma function.
We set the radius of $B_i$ to be $R_{B_i}= \max_{j \in S_i}r_{ij} + \text{SD}(\mathbf{x}_i)$ in our simulations.
The packing density considers each point in the neighborhood of $\mathbf{x}_i$ to occupy a volume of a ball of
radius equal to its separation distance. A local quality metric is defined as the relative difference
between this packing density and its optimal value,
\begin{equation}
	\text{PD}(\mathbf{x}_i) = \frac{\left|\eta_i -\eta_{opt}\right|}{\eta_{opt}},
\end{equation}
where $\eta_{opt}$ denotes an \emph{optimal packing density}, derived in sphere packing problems.
In 2D, we have $\eta_{opt} = \frac{\pi}{\sqrt{12}} \approx 0.9069$, and in 3D,
$\eta_{opt} = \frac{\pi}{\sqrt{18}}  \approx 0.7405$ \cite{conway2013sphere}.



\subsubsection{Local regularity metrics}\label{localRegularity}

The last two geometry-based quality metrics that we consider are based on distances to neighboring points.
Consider the set of all such distances
\begin{equation}
	\mathcal{D}_i	=
	\big\{ \|\mathbf{x}_i - \mathbf{x}_j\| : \mathbf{x}_j \in \mathbf{S}_i,\, j \neq i  \big\}.
\end{equation}
From this set, two local quality metrics are defined as
\begin{equation}
	\text{ALR}(\mathbf{x}_i) =\mean(\mathcal{D}_i),
	\qquad \text{RLR}(\mathbf{x}_i) =\range(\mathcal{D}_i),
\end{equation}
where ALR stands for the average of local regularity, and RLR denotes the range of local regularity.
Both these metrics have been used as point cloud quality metrics in meshfree literature with the expectation that
large discrepancies in the set of distances $\mathcal{D}_i$ indicate poor point clouds \cite{slak2019generation, van2021fast}.

\subsection{Mesh-based quality metrics}

Since mesh quality is so widely studied, it is natural to examine how mesh quality metrics translate
to indicators of point cloud quality. These types of metrics are computed by applying well-established
mesh quality metrics on a mesh generated from a point cloud. This can be done
globally by meshing the entire point cloud at once, or locally, by meshing the point cloud in each
support domain $\mathbf{S}_i$. The locally defined meshes do not impose any restriction of the local meshes
stitching together to form a global mesh~\cite{suchde2023point,suchde2017flux}.
The additional step of mesh generation can make these metrics
computationally more expensive than many other metrics used in this work.
We consider the classic aspect ratio, as the most common mesh quality metric,
a Jacobian-based metric adopted to our setting of local supports $S_i$,  along with the mass centers
deviation, which has already been used as a criterion for generating regular point clouds in the
Smoothed Particle Hydrodynamics.
These local metrics are computed on Delaunay tessellations of the point clouds.

We note that for each metric of this type, the metric value might depend only on a subset of the support
$\mathbf{S}_i$. This is in contrast with all other metric types, where the entire support is used. 

\subsubsection{Inverse aspect ratio (IAR)}\label{aspectRationAR}

The aspect ratio is one of the most widely used quality metrics in the meshing community. The normalized inverse aspect ratio
of a triangle (2D) or tetrahedron (3D) $T$ is given by
\begin{equation}
	\gamma_{T} = \kappa\frac{\rho_T}{h_T}.
	\label{aspectRatio_App}
\end{equation}
where $h_T$ is the length of the largest edge (also referred to as diameter), $\rho_T$ denotes the
inradius, and \( \kappa \) is set to \( 2\sqrt{3} \) and \( 2\sqrt{6} \) for 2D and 3D, respectively, see e.g.\
\cite{george2019meshing}.

The aspect ratio evaluates the shape of mesh elements. 
Following \cite{geuzaine2009gmsh}, we use the reciprocal of the usual aspect ratio convention.
A lower \(\gamma_T\) indicates an elongated or distorted element, while a higher \(\gamma_T\) indicates a more equilateral shape, 
which is preferred in mesh-based simulations. With the normalization in Equation~\eqref{aspectRatio_App}, valid elements satisfy 
\(0 < \gamma_T \leq 1\), where \(\gamma_T = 1\) corresponds to an equilateral triangle or regular tetrahedron. 

In order to use the aspect ratio to assess point cloud quality, we start with a triangulation of the point cloud,
which could be either locally or globally determined. 
At a point $\mathbf{x}_i$, the inverse aspect ratio (IAR) as a local point cloud quality metric is defined as
\begin{equation}\label{newAR}
	\text{IAR}(\mathbf{x}_i) = \mean_{T \in \mathbb{T}_i} \, \gamma_{T} ,
\end{equation}
where $\mathbb{T}_i$ is the set of all simplices incident on the point $\mathbf{x}_i$.

\subsubsection{Determinant of Jacobian Matrix (DJM)}\label{jacobM}

The Jacobian matrix is another widely used tool for evaluating mesh element distortion \cite{knupp2003algebraic}.
For linear elements, the Jacobian matrix defines the affine transformation from a reference element to the
physical element. For a triangle or tetrahedron $T$ with vertices $\mathbf{v}_0,\ldots,\mathbf{v}_\nu$,
the Jacobian matrix \( J_T \) is given by:
\begin{equation}
	J_T = [\mathbf{v}_{j}-\mathbf{v}_0]_{j = 1}^{\nu},
\end{equation}
where the columns of $J_T$ are edge vectors of \( T \).
While several mesh quality metrics based on $J_T$ have been considered, we restrict the discussion to the determinant
$\det(J_T)$ which indicates how much the element $T$ is distorted in comparison to the ideal reference element.

As with the IAR above, to use the Jacobian matrix to assess point cloud quality, we start with computing a triangulation of the point cloud.
At a point $\mathbf{x}_i$, we define the local point cloud quality metric as the average of the Jacobian
determinants of all elements incident on $\mathbf{x}_i$
\begin{equation}
	\text{DJM}(\mathbf{x}_i) = \mean_{T \in \mathbb{T}_i} \, \det(J_T) .
\end{equation}
%



The sign of \(\det(J_T)\) depends on the ordering of the element nodes. 
We use a consistent local numbering so that \(\det(J_T)>0\) is always satisfied.

\subsubsection{Mass centers deviation (MCD)}\label{massCenterD}

The last mesh-based metric that we consider, based on the Voronoi tessellation of the point cloud,
is derived from a point cloud generation procedure used in the Smoothed Particle Hydrodynamics (SPH) community
\cite{fu2019optimal,fu2017physics, fernandez2020delta}.
In this approach, points are iteratively repositioned so that each point coincides with the centroid (center of mass) of its
corresponding Voronoi cell, resulting in a centroidal Voronoi tessellation (CVT).

We adapt this procedure to consider the distance of a point from the centroid of the corresponding Voronoi cell.
A smaller distance would mean that the local point cloud is close to the desired CVT configuration.
For this, we define the local metric for a point for each point $\mathbf{x}_{i}$ as
%
%
\begin{equation}
\text{MCD}(\mathbf{x}_{i}) = \|{\mathbf{x}_{i} - \mathbf{z}_{i}}\|,
\end{equation}
where $\mathbf{z}_{i}$ is the centroid (mass center) of the Voronoi cell containing the point $\mathbf{x}_{i}$.


\subsection{System matrix based quality metrics}\label{system}

As the last type of metrics, we propose five metrics that use the coefficients $w_{i,j}$ of the system matrix
\eqref{discritPoisson}, in addition to the positions of the points in $\mathbf{S}_i$. These coefficients stem from the
discrete numerical differentiation operators used
in meshfree schemes. Specifically, since we test point cloud quality for meshfree collocation methods, the
metrics proposed here are based on the collocation procedures. Note that these metrics can only be computed after the system
matrix has been assembled.

\subsubsection{Growth functions (GFL1 and GFL2)}\label{growthF}

The consistency error between Equation~\eqref{PoissonEq} and its discretized form (Equation~\eqref{discritPoisson}) is defined as
follows:
\begin{equation}
	D u(\mathbf{x}_i) -  \sum_{\mathbf{x}_{j} \in \mathbf{S}_i}w_{i,j}\hat{u}_{j}.
	\label{consistError}
\end{equation}
As stated in \cite{davydov2018minimal}, the consistency error in Equation~\eqref{consistError} is bounded as follows:
\begin{equation}\label{consErEst}
	\Big| D u(\mathbf{x}_i) -  \sum_{\mathbf{x}_{j} \in \mathbf{S}_i}w_{i,j}\hat{u}_{j}\Big| \leq
	 \lVert \mathbf{w}_i \rVert_{1,p} \lVert f\lVert_{F_p}, \quad f \in F_p.
\end{equation}
Here, $p$ is the polynomial order of the method, \( F_p\) is a (Sobolev) function space with a (semi) norm \(\lVert \cdot \rVert_{F_p}\) and
\( \lVert \mathbf{w}_i \rVert_{1,p}  := \sum_{\mathbf{x}_{j} \in \mathbf{S}_i}
|w_{i,j}| \rVert \mathbf{x}_{j} - \mathbf{x}_i \lVert^{p}\) is a weighted \(\ell_1\)-norm of the weight vector
$\mathbf{w}_i=[w_{i,j}]_{\mathbf{x}_{j}\in \mathbf{S}_i}$. This motivates the local quality metric
\begin{equation}\label{GFL1}
\text{GFL1}(\mathbf{x}_{i}) = \|\mathbf{w}_i\|_{1,p}.
\end{equation}
Since $\|\mathbf{w}_i\|_{1,p}\le \sqrt{|\mathbf{S}_{i}|}\,\|\mathbf{w}_i\|_{2,p}$ for
\( \lVert \mathbf{w}_i \rVert_{2,p}^2  := \sum_{\mathbf{x}_{j} \in \mathbf{S}_i}
w_{i,j}^2 \rVert \mathbf{x}_{j} - \mathbf{x}_i \lVert^{2p}\)
and the supports $\mathbf{S}_{i}$ are
relatively small, we also consider the metric
\begin{equation}\label{GFL2}
\text{GFL2}(\mathbf{x}_{i}) = \|\mathbf{w}_i\|_{2,p}.
\end{equation}
In the case of the Poisson problem we minimize $\|\mathbf{w}_i\|_{2,p}$ subject to the condition that the
numerical differentiation formula \eqref{diffFormula} satisfies the polynomial exactness of order \(p\), see Section
\ref{numericalMethods}. As shown in \cite{davydov2018minimal}, this minimal value of $\|\mathbf{w}_i\|_{2,p}$
represents a 2-norm \emph{growth function}, which gives the names for the two metrics. Remarkably, the growth function also
appears in the estimates of the numerical differentiation by the kernel-based methods \cite{davydov2016error} and therefore
can be considered as a measure of the quality of the local point set $\mathbf{S}_{i}$ with respect to its ability to provide
accurate numerical differentiation of the operator $D$ with polynomial order $p$.

\subsubsection{Condition number of the weighted Vandermonde matrix (CN)}\label{conditionN}

Consider a point $\mathbf{x}_i$ with support  $ \mathbf{S}_i$ and a space of
multivariate polynomials of order \(p\), where the dimension of this polynomial space is \(P =\binom{p - 1 + \nu}{\nu}\).
The weighted Vandermonde matrix  \(V_i\) for the point $\mathbf{x}_i$ is an \(m \times P\) matrix, with $m=|\mathbf{S}_i|$,
that contains the values of the monomials of degree up to $p-1$ at the points $\mathbf{x}_{j}\in \mathbf{S}_i$.
This matrix appears in the formulation of the least squares numerical differentiation problem that we use to determine the
entries $w_{i,j}$ of the system matrix.

The sensitivity of the least squares solution depends on the condition number of \(V_i\). Moreover,
the size of this condition number is an important criterion for assessing the numerical stability of the
method. The \(\ell_2\)-condition number of \(V\) is defined as:
\begin{equation}
	\text{CN}(\mathbf{x}_i) =  \frac{\sigma_{\max}(V_i)}{\sigma_{\min}(V_i)},
\end{equation}
where \(\sigma_{\max}(V_i)\) and \(\sigma_{\min}(V_i)\) are the largest and the smallest non-zero singular
values of \(V\), respectively. In this work, we use the condition number CN as a quality metric
for the point cloud.

\subsubsection{Recovery error (RE)}\label{maxAppErro}

Let \( K : \mathbb{R}^\nu \times \mathbb{R}^\nu \to \mathbb{R} \) be a conditionally positive
definite kernel of order \( p \), associated with a native space \( \mathcal{F}_{K,p} \).
It has been shown recently in \cite{davydov2021approximation}  that for any numerical differentiation
formula \eqref{diffFormula} exact for polynomials of order \( p \),  the
worst-case recovery error for functions in \( \mathcal{F}_{K,p} \), defined by
\begin{equation}
	\tau_i := \sup_{\substack{f \in \mathcal{F}_{K,p} \\ \lVert f
	\lVert_{K,p} \leqslant 1}} \Big| D f(\mathbf{x}_i) - \sum_{\mathbf{x}_{j} \in
	\mathbf{S}_i}
	w_{i,j} f(\mathbf{x}_{j}) \Big|,
\end{equation}
can be computed by
\begin{equation}\label{Qbf}
		\tau^{2}_i = D'D''K(\mathbf{x}_{i},\mathbf{x}_{i})
		- \sum_{\mathbf{x}_{j} \in \mathbf{S}_i}
		w_{i,j} \left( D'K(\mathbf{x}_{i},\mathbf{x}_{j}) +
D''K(\mathbf{x}_{j},\mathbf{x}_{i}) \right)  + \sum_{\mathbf{x}_{j}, \mathbf{x}_{k} \in \mathbf{S}_i} w_{i,j} w_{i,k}
		K(\mathbf{x}_{j},\mathbf{x}_{k}),
\end{equation}
where \( D' \) and \( D'' \) denote the differentiation operator $D$ acting on the first and the second arguments
of the kernel \( K\), respectively. Note that \eqref{Qbf} has been known before \cite{davydov2021approximation}
under certain additional conditions on the support $\mathbf{S}_i$, see e.g.~\cite{wendland2004scattered,davydov2016error}.
Similar to the growth functions of Section~\ref{growthF}, small $\tau_i$ indicates the ability of the local point set $\mathbf{S}_i$ to provide
small consistency error \eqref{consistError}.

In this study we choose the \emph{polyharmonic kernel} $K(\mathbf{x},\mathbf{y})=\|\mathbf{x}-\mathbf{y}\|^{2p-1}$,
and define the point quality metric
$$
\text{RE}(\mathbf{x}_i)=\tau_i,$$

\subsubsection{Stability-based quality metric (ST)}\label{stabilityMe}
We propose a quality metric based on the stability constant of a local square submatrix $A_i$ of the system matrix
$A=[w_{j,k}]_{\mathbf{x}_j,\mathbf{x}_k\in\Omega}$, obtained from the rows and columns of $A$ corresponding to all
$\mathbf{x}_k\in\Omega$ such that $\mathbf{x}_k\in \mathbf{S}_j$ for some $\mathbf{x}_j\in \mathbf{S}_i\cap \Omega$.
In other words, we collect interior points $\mathbf{x}_k$ in the union of the supports $\mathbf{S}_j$ of all interior neighboring points
$\mathbf{x}_j\in \mathbf{S}_i$, and build $A_i$ from the corresponding rows and columns of $A$.

The matrix $A_i$ is the system matrix of a local homogeneous Dirichlet problem in a domain with the interior corresponding to
the  above points $\mathbf{x}_k$, and homogeneous boundary conditions enforced on all other points in their supports
$\mathbf{S}_k$.


The local stability-based quality metric is defined as
\begin{equation}
	\text{ST}(\mathbf{x}_i) = \| A_i^{-1} \|_{\infty}.
	\label{stabilityConstant}
\end{equation}
We estimate the stability constant by using the formula
$\| A_i^{-1} \|_{\infty}=\cond_1(A_i^{T})/\| A_i \|_{\infty}$ suggested in \cite{schaback2017error},
where the 1-norm condition number $\cond_1$ of a matrix is effectively estimated by
an efficient algorithm  available in MATLAB via the command \texttt{condest}.

The stability constant $\text{ST}(\mathbf{x}_i)$ multiplied by an estimate of the consistency errors like the one in
\eqref{consErEst} provides a rigorous error bound for the meshfree solution of the homogeneous Dirichlet problem. Therefore,
a smaller ST may be seen as an indicator of a higher quality of a point cloud.

\subsubsection{Diagonal weights (DW)}\label{sysMat}
It is known, see e.g.~\cite{hackbusch2017elliptic} that a big advantage for a classical finite
difference discretization of a Poisson equation is a diagonally dominant system matrix. Also in the meshfree setting
the central weight value \(w_{i,i}\) should be significantly
larger than the absolute values of the other entries \(w_{i,j}\) for $\mathbf{x}_j\in \mathbf{S}_i$, \(j \neq i\)
\cite{suchde2018conservation}.
Therefore, we consider the local quality metric that just picks the diagonal weight of the system matrix
\begin{equation}
	\text{DW}(\mathbf{x}_i) = w_{i,i},
	\label{DW}
\end{equation}
and expect that a higher value of an aggregation map of it indicates a better point cloud for the Poisson equation.


\section{Additional results}
\label{sec:AppResults}

This appendix section presents additional results from the numerical studies discussed in
Section~\ref{eval}, using alternative analytical functions to those used in
the main text. See Table~\ref{testsFunction} for the different analytical functions considered.

For the parameter study on support size variation (Section~\ref{sec:Support}), the Spearman rank
correlation between the metric value and numerical error for analytical function $u_2$ is presented
in Table~\ref{tab:spearma-subsetsize-u2}. Similarly, for the polynomial order variation study
(Section~\ref{sec:POrder}), the corresponding results for analytical function $u_2$ are tabulated in
Table~\ref{tab:spearma-porder-u2}.
Furthermore, for the $\mathbb{R}^3$ Poisson problem (Section~\ref{sec:3D}), the corresponding results
for the analytical function $u_4$ are listed in Table~\ref{tab:spearma-3d-u2}.
As noted in Section~\ref{eval}, most metrics exhibit very similar performance across the
different analytical functions considered. The tables included here are provided for completeness.

\begin{table}[!ht]
	\centering
	\caption{Parameter study: varying support size with analytical solution \( u_2 \) (see
		Table~\ref{testsFunction}).
		Spearman's rank correlation values \(\rho\) between each quality metric and numerical error for
		different support sizes. The support size is given by $\varepsilon$, which represents
		the scaling factor (see Equation~\eqref{range}).  The  columns denoted \(\bar{\rho}\) and $\bar\sigma$ give the
	averages of $\rho$ and of the performance score $\sigma$ across all four scaling factors. The other parameters are set to
		polynomial order \( p = 3 \) and $n = 1600$ points in the domain.} 
	\label{tab:spearma-subsetsize-u2}
	\begin{tabular}{lcccccccc}
		\toprule
		\textbf{Metric} & \textbf{Map} & \textbf{Type} & \multicolumn{4}{c}{\(\rho\)}& \textbf{\( \bar{\rho}\)} & $\bar\sigma$\\
		\cmidrule(lr){4-7}
		& & & \textbf{$\varepsilon =1.8$}&\textbf{$\varepsilon =2.1$}&\textbf{$\varepsilon = 2.5$}
		&\textbf{$\varepsilon=2.8$} & \\
		\midrule
	RE   & std    & system      & 0.921 & 0.859 & 0.860 & 0.864 & 0.876 & 1 \\
	IMPF  & rms    & equilibrium & 0.919 & 0.859 & 0.861 & 0.865 & 0.876 & 1 \\
	GFL1  & std    & system      & 0.920 & 0.860 & 0.860 & 0.864 & 0.876 & 1 \\
	DJM   & std    & mesh        & 0.919 & 0.858 & 0.860 & 0.863 & 0.875 & 2 \\
	GFL2  & std    & system      & 0.920 & 0.857 & 0.859 & 0.862 & 0.875 & 2 \\
	FD    & std    & geometry        & 0.919 & 0.857 & 0.858 & 0.861 & 0.874 & 3 \\
	DW    & mean   & system      & 0.918 & 0.857 & 0.858 & 0.862 & 0.874 & 3 \\
	LKPF   & median & equilibrium & 0.918 & 0.853 & 0.860 & 0.861 & 0.873 & 4 \\
	ALR   & std    & geometry    & 0.916 & 0.856 & 0.861 & 0.860 & 0.873 & 4 \\
	WKPF  & mean   & equilibrium & 0.919 & 0.856 & 0.858 & 0.859 & 0.873 & 4 \\
	RKPF  & median & equilibrium & 0.914 & 0.854 & 0.859 & 0.861 & 0.872 & 5 \\
	CE    & std    & energy      & 0.914 & 0.854 & 0.860 & 0.859 & 0.872 & 5 \\
	RKE  & median & energy      & 0.915 & 0.854 & 0.859 & 0.860 & 0.872 & 5 \\
	IME   & median & energy      & 0.915 & 0.854 & 0.858 & 0.859 & 0.872 & 5 \\
	MCD   & std    & mesh        & 0.915 & 0.854 & 0.858 & 0.858 & 0.871 & 6 \\
	RLR   & mean   & geometry    & 0.912 & 0.851 & 0.859 & 0.859 & 0.870 & 7 \\
	CPF   & mean   & equilibrium & 0.919 & 0.844 & 0.856 & 0.859 & 0.869 & 8 \\
	WKE   & std    & energy      & 0.917 & 0.853 & 0.850 & 0.846 & 0.867 & 10 \\
	LKE   & std    & energy      & 0.914 & 0.852 & 0.849 & 0.847 & 0.865 & 12 \\
	IAR    & std    & mesh        & 0.897 & 0.840 & 0.849 & 0.847 & 0.858 & 19 \\
	CN    & median & system      & 0.886 & 0.845 & 0.850 & 0.849 & 0.857 & 20 \\
	SD    & rms    & geometry    & 0.896 & 0.839 & 0.848 & 0.843 & 0.857 & 20 \\
	PD    & mean   & geometry    & 0.895 & 0.839 & 0.847 & 0.843 & 0.856 & 21 \\
	UR    & median & geometry        & 0.893 & 0.833 & 0.839 & 0.835 & 0.850 & 27 \\
	ST    & max    & system      & 0.856 & 0.793 & 0.799 & 0.808 & 0.814 & 63 \\
		\bottomrule
	\end{tabular}
\end{table}

\begin{table}[!ht]
	\centering
	\caption{Parameter study: varying polynomial order with analytical solution \( u_2 \) (see
		Table~\ref{testsFunction}).
		Spearman's rank correlation values \(\rho\) between each quality metric and numerical error for
		different polynomial orders. The polynomial order is given by $p$.  The  columns denoted \(\bar{\rho}\) and $\bar\sigma$ give the
	averages of $\rho$ and of the performance score $\sigma$ across all polynomial orders. The other parameters are set to scaling factors of 2.1, 2.3, and 2.8 
	for polynomial orders $p = 3, 4,$ and \(5\), respectively, with $n = 1600$ points in the domain.} 
	\label{tab:spearma-porder-u2}
	\begin{tabular}{lccccccc}
		\toprule
		\textbf{Metric} & \textbf{Map} & \textbf{Type} & \multicolumn{3}{c}{\(\rho\)}&\textbf{\( \bar{\rho}\)} & $\bar\sigma$\\
		\cmidrule(lr){4-6}
		& & & \textbf{$p =3$} & \textbf{$p=4$} &\textbf{$p =5$} & \\
		\midrule
	RE   & std    & system      & 0.859 & 0.972 & 0.951 & 0.927 & 1 \\
	DJM   & std    & mesh        & 0.858 & 0.971 & 0.948 & 0.926 & 2 \\
	IMPF  & rms    & equilibrium & 0.859 & 0.970 & 0.947 & 0.925 & 3 \\
	GFL2  & mean   & system      & 0.857 & 0.969 & 0.950 & 0.925 & 3 \\
	FD    & rms    & geometry        & 0.857 & 0.970 & 0.948 & 0.925 & 3 \\
	RKPF  & median & equilibrium & 0.854 & 0.971 & 0.946 & 0.924 & 4 \\
	GFL1  & mean   & system      & 0.850 & 0.970 & 0.949 & 0.923 & 5 \\
	RKE  & median & energy      & 0.854 & 0.970 & 0.946 & 0.923 & 5 \\
	IME   & median & energy      & 0.854 & 0.969 & 0.945 & 0.923 & 5 \\
	LKPF   & mean   & equilibrium & 0.855 & 0.968 & 0.945 & 0.922 & 6 \\
	WKPF  & mean   & equilibrium & 0.856 & 0.967 & 0.944 & 0.922 & 6 \\
	MCD   & std    & mesh        & 0.854 & 0.967 & 0.944 & 0.921 & 7 \\
	RLR   & median & geometry    & 0.851 & 0.967 & 0.943 & 0.920 & 8 \\
	ALR   & std    & geometry    & 0.856 & 0.965 & 0.939 & 0.920 & 8 \\
	CE    & std    & energy      & 0.854 & 0.965 & 0.939 & 0.920 & 8 \\
	CPF   & mean   & equilibrium & 0.844 & 0.968 & 0.945 & 0.919 & 9 \\
	WKE   & std    & energy      & 0.853 & 0.959 & 0.936 & 0.916 & 12 \\
	LKE   & std    & energy      & 0.852 & 0.954 & 0.934 & 0.913 & 15 \\
	IAR    & std    & mesh        & 0.840 & 0.951 & 0.932 & 0.908 & 20 \\
	SD    & rms    & geometry    & 0.840 & 0.949 & 0.922 & 0.904 & 24 \\
	CN    & median & system      & 0.845 & 0.936 & 0.930 & 0.904 & 24 \\
	PD    & mean   & geometry    & 0.839 & 0.949 & 0.921 & 0.903 & 25 \\
	UR    & median & geometry        & 0.833 & 0.943 & 0.915 & 0.897 & 31 \\
	DW    & std    & system      & 0.822 & 0.941 & 0.923 & 0.895 & 33 \\
	ST    & median & system      & 0.849 & 0.741 & 0.950 & 0.847 & 81 \\
		\bottomrule
	\end{tabular}
\end{table}

\begin{table}[!ht]
	\centering
	\caption{Poisson problem in \(\mathbb{R}^3\) with analytical solution \( u_4 \) (see
		Table~\ref{testsFunction}):
		The Spearman's rank correlation \(\rho\) quantifying the monotonic dependence between quality metrics
		and numerical error. The discretization resolution is given by \(n\), representing the
		number of points in the point clouds considered.  The  columns denoted \(\bar{\rho}\) and $\bar\sigma$ give the
	averages of $\rho$ and of the performance score $\sigma$ across both resolutions. Other parameters are set to polynomial order \( p = 3 \) and
		scaling factor \( \varepsilon = 2.1 \).} 
	\label{tab:spearma-3d-u2}
	\begin{tabular}{lcccccc}
		\toprule
		\textbf{Metric} & \textbf{Map} & \textbf{Type} & \multicolumn{2}{c}{\(\rho\)} & \textbf{\(\bar{\rho}\)} & $\bar\sigma$\\
		\cmidrule(lr){4-5}
		& & & \( n = 2744 \) & \( n = 39304 \) & \\
		\midrule
	RE   & std    & system      & 0.944 & 0.982 & 0.963 & 1 \\
	GFL2  & std    & system      & 0.942 & 0.982 & 0.962 & 2 \\
	GFL1  & std    & system      & 0.942 & 0.982 & 0.962 & 2 \\
	FD    & std    & geometry        & 0.942 & 0.982 & 0.962 & 2 \\
	DJM   & std    & mesh        & 0.939 & 0.982 & 0.960 & 4 \\
	IMPF  & rms    & equilibrium & 0.938 & 0.982 & 0.960 & 4 \\
	DW    & mean   & system      & 0.937 & 0.982 & 0.960 & 4 \\
	RKPF  & median & equilibrium & 0.937 & 0.979 & 0.958 & 6 \\
	CN    & median & system      & 0.934 & 0.981 & 0.957 & 7 \\
	IAR    & std    & mesh        & 0.933 & 0.981 & 0.957 & 7 \\
	IME   & median & energy      & 0.932 & 0.981 & 0.957 & 7 \\
	CPF   & mean   & equilibrium & 0.935 & 0.974 & 0.954 & 10 \\
	RKE  & median & energy      & 0.929 & 0.977 & 0.953 & 11 \\
	RLR   & mean   & geometry    & 0.923 & 0.976 & 0.949 & 15 \\
	LKPF   & mean   & equilibrium & 0.912 & 0.982 & 0.947 & 17 \\
	WKPF  & mean   & equilibrium & 0.911 & 0.982 & 0.947 & 17 \\
	WKE   & std    & energy      & 0.921 & 0.970 & 0.945 & 19 \\
	MCD   & std    & mesh        & 0.919 & 0.969 & 0.944 & 20 \\
	LKE   & std    & energy      & 0.900 & 0.975 & 0.937 & 27 \\
	UR    & range  & geometry        & 0.921 & 0.951 & 0.936 & 28 \\
	SD    & max    & geometry    & 0.893 & 0.941 & 0.917 & 47 \\
	ALR   & std    & geometry    & 0.850 & 0.981 & 0.916 & 48 \\
	CE    & std    & energy      & 0.848 & 0.981 & 0.915 & 49 \\
	PD    & std    & geometry    & 0.891 & 0.934 & 0.912 & 52 \\
	ST    & range  & system      & 0.847 & 0.897 & 0.872 & 92 \\
		\bottomrule
	\end{tabular}
\end{table}
\end{appendices}

	\section*{Acknowledgements}
	Pratik Suchde would like to acknowledge funding from the Institute of Advanced Studies,
	University of Luxembourg, under the AUDACITY programme.

	\setlength{\bibsep}{0pt plus 0.3ex}
	\bibliographystyle{abbrv} 
 	\bibliography{mybib}{}
 \end{document}